\documentclass[12pt]{article} 
\usepackage{t1enc}
\usepackage[latin1]{inputenc}
\usepackage[english]{babel}
\usepackage{amsmath,amsthm}
\usepackage{amsfonts}
\usepackage{latexsym}
\usepackage[dvips]{graphicx}
\usepackage{graphicx}
\usepackage[natural]{xcolor}
\usepackage{epstopdf}
\theoremstyle{plain}
\newtheorem{thm}{Theorem}[section]

\newtheorem{claim}[thm]{Claim}
\newtheorem{theorem}{Theorem}
\numberwithin{theorem}{section}
\newtheorem{lemma}{Lemma}
\numberwithin{lemma}{section}

\numberwithin{corollary}{section}

\numberwithin{conjecture}{section}
\newtheorem*{conjecture*}{Conjecture}

\title{A note on the edge choosability of $K_{5}$-minor free graphs}
\author{Jieru Feng, Jianliang Wu\thanks{Corresponding author. \emph{E-mail address:} jlwu@sdu.edu.cn}, Fan Yang\\
{\small {School of Mathematics, Shandong University, Jinan, 250100, China}}}
\date{}
\begin{document}
\baselineskip 0.7cm

\maketitle
\begin{abstract}
For a planar graph $G$, Borodin stated that $G$ is $(\Delta+1)$-edge-choosable if $\Delta \geq9$ and later Bonamy showed that $G$ is $9$-edge-choosable if $\Delta=8$. At the same time, Borodin et al. proved that $G$ is $\Delta$-edge-choosable if $\Delta\geq12$. In the paper, we extend these results to  $K_5$-minor free graphs.

\end{abstract}

\section{Introduction}

Throughout the paper, the graph $G$ we considered is finite, undirected and simple, all undefined notation and definitions can be found in \cite{B}. Let $V(G)$ (resp. $E(G)$) be the vertex set (edge set) of $G$. For a vertex $v\in V(G)$, let $N_G(v)$ (or simply $N(v)$) denote the set of vertices adjacent to $v$ and $d(v) = |N(v)|$ denote the degree of $v$. Moreover, we use $\Delta(G)$ and $\delta(G)$ (or simply $\Delta$ and $\delta$) to denote the maximum degree and the minimum degree of $G$, respectively. A vertex of degree $k$ (at least $k$, or at most $k$) is called a $k$-vertex ($k^{+}$-vertex, or $k^{-}$-vertex). A $k$-$cycle$ is a cycle of length $k$. In particular, a $3$-cycle is also said to be a \emph{triangle}.

A \emph{$k$-edge coloring} of a graph $G=(V, E)$ is a mapping $\phi:E(G)\rightarrow \{1, 2, \ldots, k\}$ such that every two adjacent edges receive different colors. The \emph{edge chromatic number} $\chi'(G)$ of $G$ is the smallest integer $k$ such that $G$ has a $k$-edge coloring. In 1964, Vizing \cite{V} showed that $\Delta \leq \chi'(G) \leq \Delta+1$ admits for every graph $G$. Specially, he also proved every planar graph $G$ with $\Delta(G)\geq8$ satisfies $\Delta(G)=\chi'(G)$  \cite{V21}. In addition, every planar graph $G$ with $\Delta(G)=7$ satisfies $\Delta(G)=\chi'(G)$, which was proved by Sanders and Zhao \cite{S}, independently by Zhang \cite{zhang}.

For any $L : E(G)\rightarrow \mathcal{P}(N)$ list assignment of colors to the edges of a graph $G$, the graph $G$ is \emph{L-edge-colorable} if there exists an proper edge coloring of $G$ such that the color of every edge $e \in E(G)$ belongs to $L(e)$. A graph $G$ is said to be \emph{k-edge-choosable} (or \emph{list k-edge-colorable}) if $G$ is $L$-edge-colorable for any list assignment $L$ such that $|L(e)|\geq k$ for each edge $e\in E(G)$. We denote by $\chi_{\ell}'(G)$ the smallest $k$ such that $G$ is $k$-edge-choosable.

Note that if all the lists are equal, then the edge coloring is a special case of list edge coloring. Thus $\chi'(G)\leq\chi_{\ell}'(G)$. In \cite{Jensen}, Jensen et al. conjectured that $\chi'(G)=\chi_{\ell}'(G)$ and the conjecture remains wide open. However, for some special classes of graphs the conjecture is true, for example, planar graphs of the maximum degree at least 12. In \cite{Borodin2}, Borodin et al. stated that $\chi_{\ell}'(G)=\Delta$ for every planar graph $G$ with $\Delta(G)\geq 12$.

By $\Delta\leq\chi'(G)\leq\Delta(G)+1$, the conjecture can be weakened into $\chi_{\ell}'(G)\leq\Delta(G)+1$ (see \cite{V1}). It is reasonable to consider the problem for planar graphs. In \cite{Borodin}, Borodin proved that every planar graph $G$ with $\Delta(G)\geq9$ satisfies $\chi_{\ell}'(G)\leq\Delta(G)+1$. More than twenty years later, Bonamy improved this result, which states that every planar graph $G$ with $\Delta(G)\leq8$ satisfies $\chi_{\ell}'(G)\leq9$ in \cite{Bonamy}. So we get that $\chi_{\ell}'(G)\leq\Delta(G)+1$ for every planar graph $G$ with $\Delta(G)\geq8$.

In the paper we consider the list edge coloring of $K_{5}$-minor free graphs. By \emph{identifying} nonadjacent vertices $x$ and $y$ of a graph $G$, we means that replacing these vertices by a single vertex incident with all the edges which were incident in $G$ with either $x$ or $y$. Let $e = xy$ be an edge of a graph $G = (V,E)$. To \emph{contract} an edge $e$ of a graph $G$ is to delete the edge first and then identify its end-vertices and finally deleting all multiple edges. A graph $H$ is a \emph{minor} of a graph $G$ if $G$ has a subgraph contractible to $H$; $G$ is called $H$-\emph{minor} \emph{free} if $G$ does not have $H$ as a minor. It is well-known that every planar graph contains neither $K_{5}$-minor nor $K_{3,3}$-minor.

Recently, the result of edge chromatic number in planar graphs has been extended to that of $K_5$-minor free graphs, which states that
if $G$ is a $K_5$-minor free graph with maximum degree $\Delta(G) \geq 7$, then $\chi'(G)=\Delta(G)$, see \cite{feng}. It is an interesting problem to determine whether the results of edge choosability in planar graphs could be extended to that of $K_5$-minor free graphs. Here we get the following theorems.

\begin{theorem}\label{th1}
Let G be a $K_{5}$-minor free graph. If $\Delta(G) \geq 8$, then $\chi'_{\ell}(G) \leq \Delta(G)+1$.
\end{theorem}

\begin{theorem}\label{th2}
Let G be a $K_{5}$-minor free graph. If $\Delta(G) \geq 12$, then $\chi'_{\ell}(G) = \Delta(G)$.
\end{theorem}


\section{Terminology and notation}
Before proceeding, we introduce the following notation. Let $G$ be a planar graph having a plane drawing and let $F$ be the face set of $G$. For a face $f$ of $G$, the $degree$ $d(f)$ is the number of edges incident with it, where each cut-edge is counted twice. A \emph{$k$-face}, \emph{$k^{-}$-face} and \emph{$k^{+}$-face} is a face of degree $k$, at most $k$ and at leat $k$, respectively. Note that in all figures of the paper, all vertices marked with $\bullet$ have no edge of $G$ incident with it other than those shown.

\begin{figure}[htbp]
\begin{center}
\includegraphics[scale=0.2]{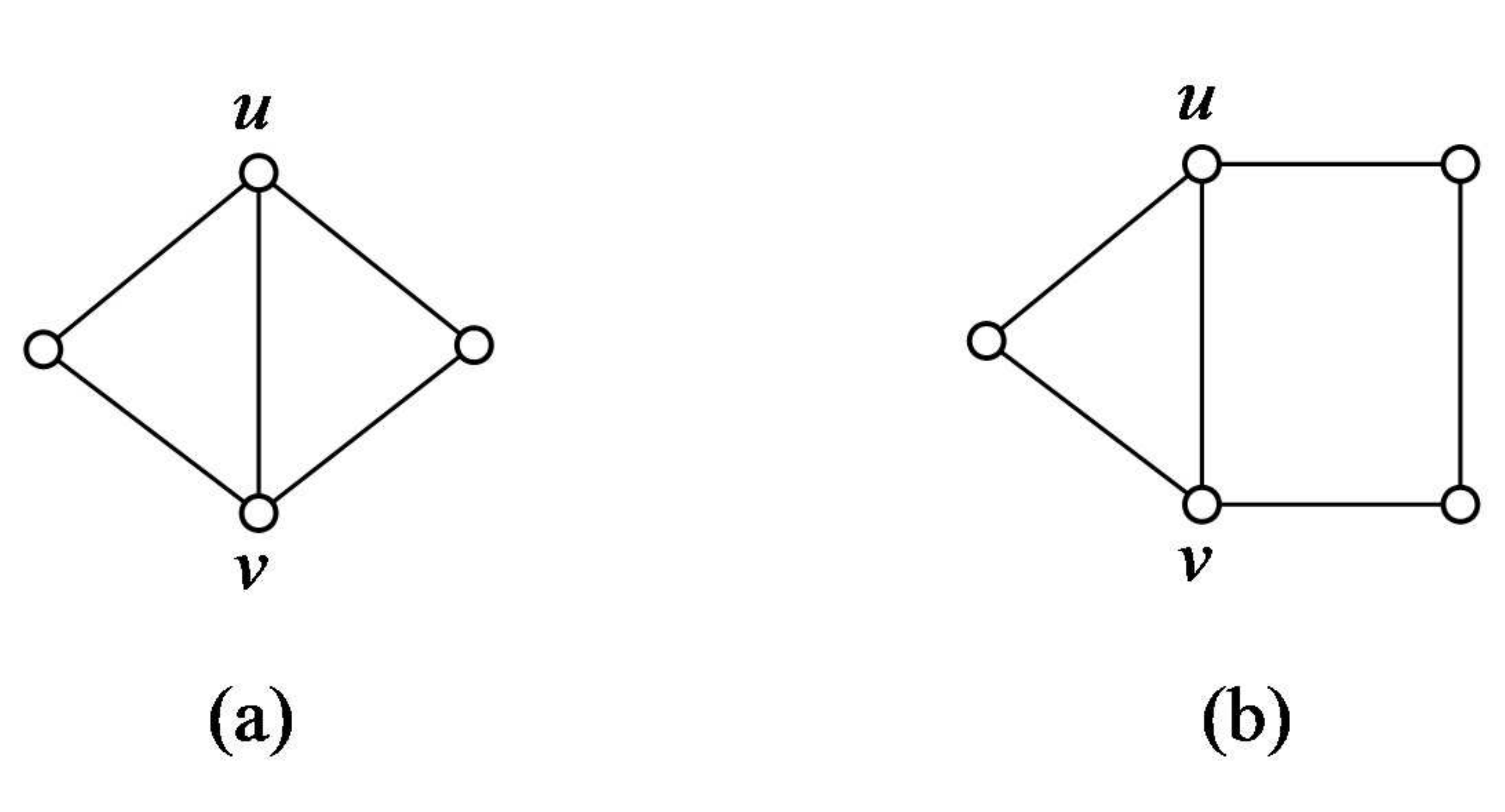}
\caption{The vertex $v$ is a weak neighbor of $u$ (a) or a semiweak neighbor of $u$ (b).}
\label{weak}
\end{center}
\end{figure}

For a given planar embedding and two vertices $u$, $v$ of $G$, $v$ is a \emph{weak} \emph{neighbor} of $u$ when the two faces incident with the edge $uv$ are triangles (see Figure \ref{weak}(a)); $v$ is a \emph{semiweak} \emph{neighbor} of $u$ when one of the two faces incident with the edge $uv$ is a triangle and the other is a 4-face (see Figure \ref{weak}(b)).

Next we consider a weak neighbor $v$ of $u$ with $d(v) = 5$ , we define some special types.

\begin{figure}[htbp]
\begin{center}
\includegraphics[scale=0.25]{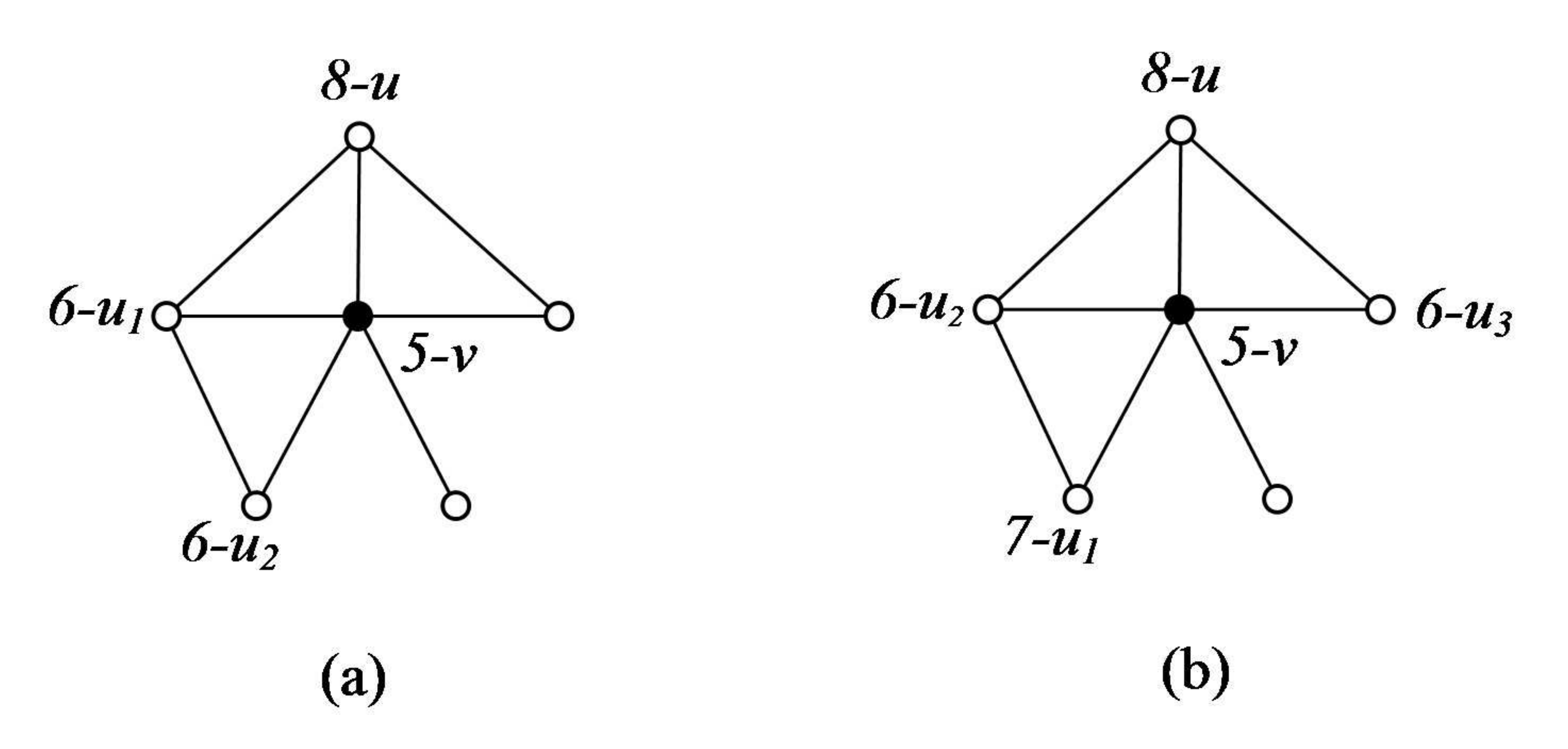}
\caption{The vertex $v$ is an $E_2$-neighbor of $u$.}
\label{E2}
\end{center}
\end{figure}

$\bullet$ The vertex $v$ is an $E_{2}$-$neighbor$ of $u$ with $d(u) = 8$ if one of the two following conditions is verified:

$\cdot$ $v$ is adjacent to two 6-vertices $u_1$ and $u_2$ such that $(u, v, u_1)$ and $(u_1, v, u_2)$ are 3-faces (see Figure \ref{E2} (a)).

$\cdot$ $v$ is adjacent to a 7-vertex $u_1$ and two 6-vertices $u_2$ and $u_3$ such that $(u_1, v, u_2)$, $(u_2, v, u)$ and $(u, v, u_3)$ are 3-faces (see Figure \ref{E2} (b)).

\begin{figure}[htbp]
\begin{center}
\includegraphics[scale=0.25]{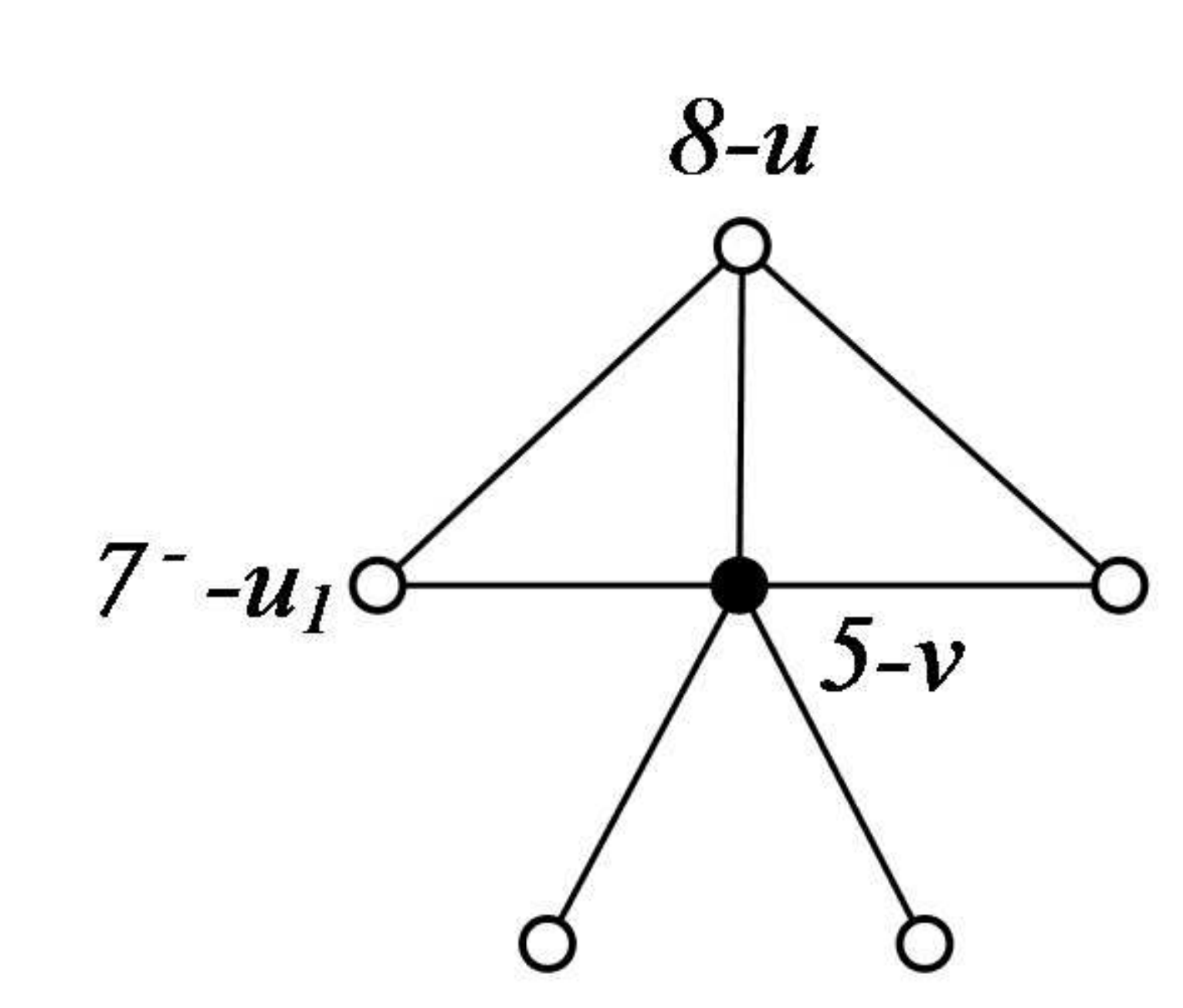}
\caption{The vertex $v$ is an $E_3$-neighbor of $u$.}
\label{E3}
\end{center}
\end{figure}

$\bullet$ The vertex $v$ is an $E_{3}$-$neighbor$ of $u$ with $d(u) = 8$ if $v$ is not an $E_2$-neighbor of $u$, and $v$ is adjacent to a $7^-$-vertex $u_1$ such that $(u, v, u_1)$ is a 3-face (see Figure \ref{E3}) .

\begin{figure}[htbp]
\begin{center}
\includegraphics[scale=0.25]{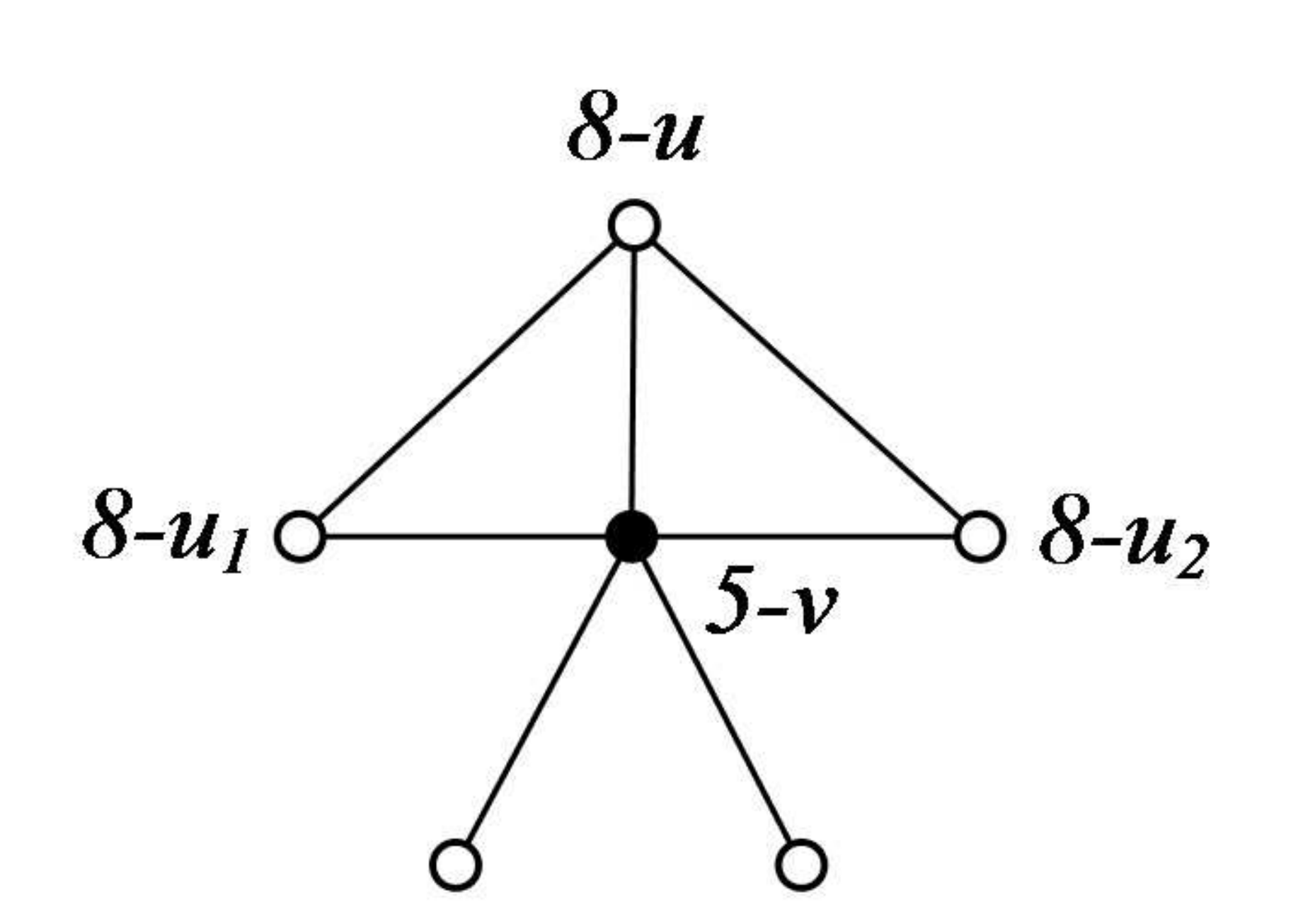}
\caption{The vertex $v$ is an $E_4$-neighbor of $u$.}
\label{E4}
\end{center}
\end{figure}

$\bullet$ The vertex $v$ is an $E_{4}$-$neighbor$ of $u$ with $d(u) = 8$ if $v$ is neither an $E_2$ nor $E_3$-neighbor of $u$, that is $v$ is adjacent to two $8$-vertex $u_1$ and $u_2$ such that $(u, v, u_1)$ and $(u, v, u_2)$ are 3-faces (see Figure \ref{E4}) .

\begin{figure}[htbp]
\begin{center}
\includegraphics[scale=0.25]{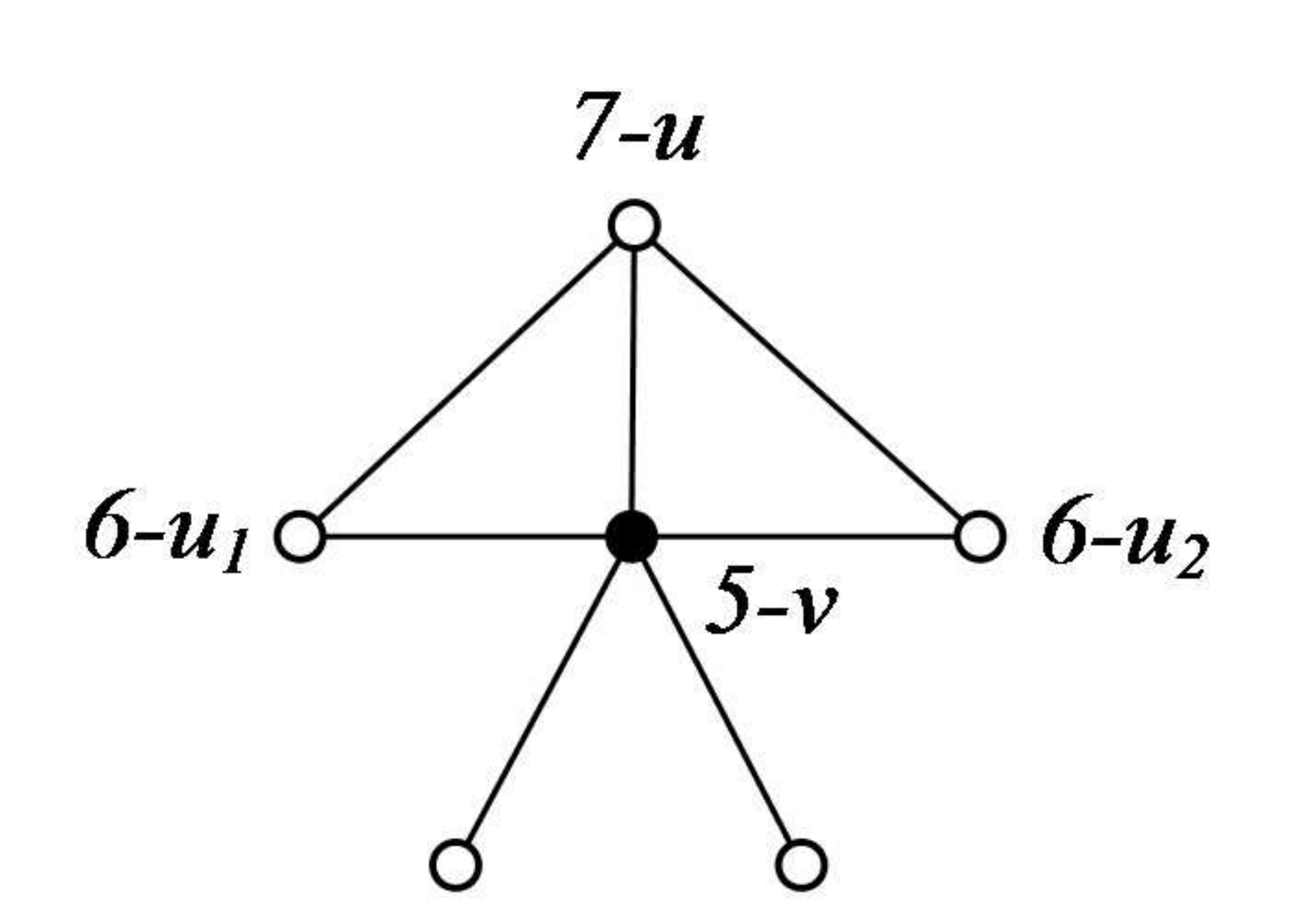}
\caption{The vertex $v$ is an $S_2$-neighbor of $u$.}
\label{S2}
\end{center}
\end{figure}

$\bullet$ The vertex $v$ is an $S_{2}$-$neighbor$ of $u$ with $d(u) = 7$ if $v$ adjacent to two $6$-vertex $u_1$ and $u_2$ such that $(u, v, u_1)$ and $(u, v, u_2)$ are 3-faces (see Figure \ref{S2}) .

\begin{figure}[htbp]
\begin{center}
\includegraphics[scale=0.25]{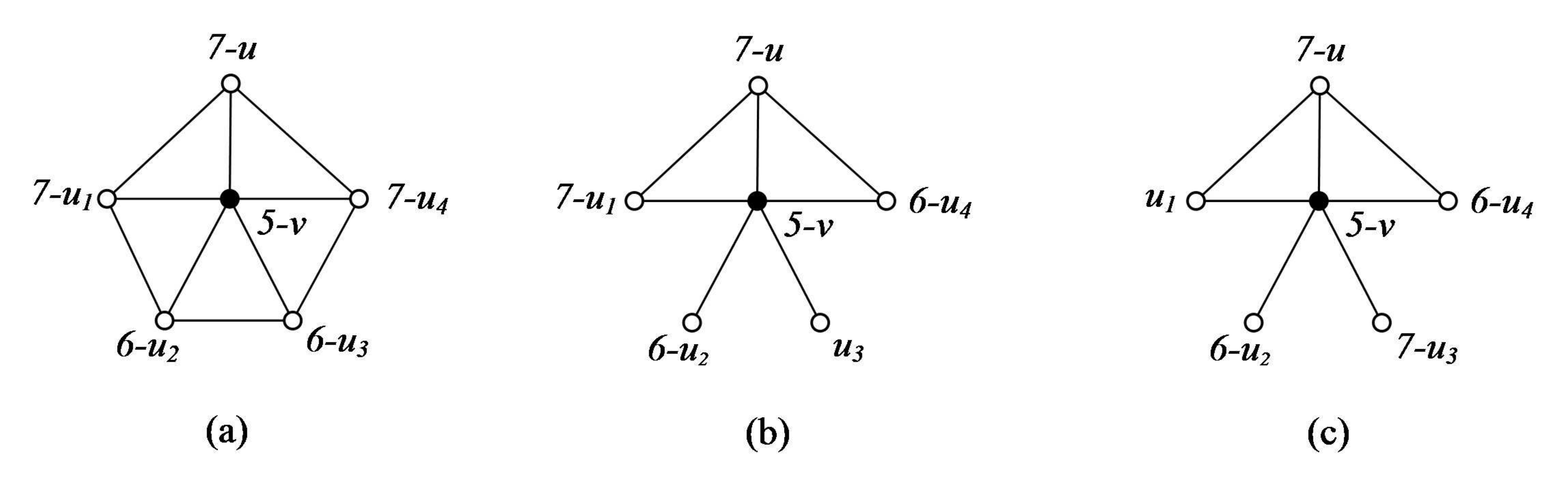}
\caption{The vertex $v$ is an $S_3$-neighbor of $u$.}
\label{S3}
\end{center}
\end{figure}

$\bullet$ The vertex $v$ is an $S_{3}$-$neighbor$ of $u$ with $d(u) = 7$ if $v$ is not an $S_2$-neighbor of $u$, and $v$ is adjacent to four vertices $u_1$, $u_2$, $u_3$ and $u_4$ such that $(u, v, u_1)$ and $(u, v, u_4)$ are 3-faces, and one of the two following conditions is verified:

$\cdot$ $d(u_1)=d(u_4)=7$, $d(u_2)=d(u_3)=6$ and $(u_1, v, u_2)$, $(u_2, v, u_3)$ and $(u_3, v, u_4)$ are 3-faces (see Figure \ref{S3} (a)).

$\cdot$ $d(u_2)=d(u_4)=6$. $d(u_1)=7$ (see Figure \ref{S3} (b)) or $d(u_3)=7$ (see Figure \ref{S3} (c)). Note that there is no constraint on the order of $u_2$ and $u_3$ in the embedding.

\begin{figure}[htbp]
\begin{center}
\includegraphics[scale=0.25]{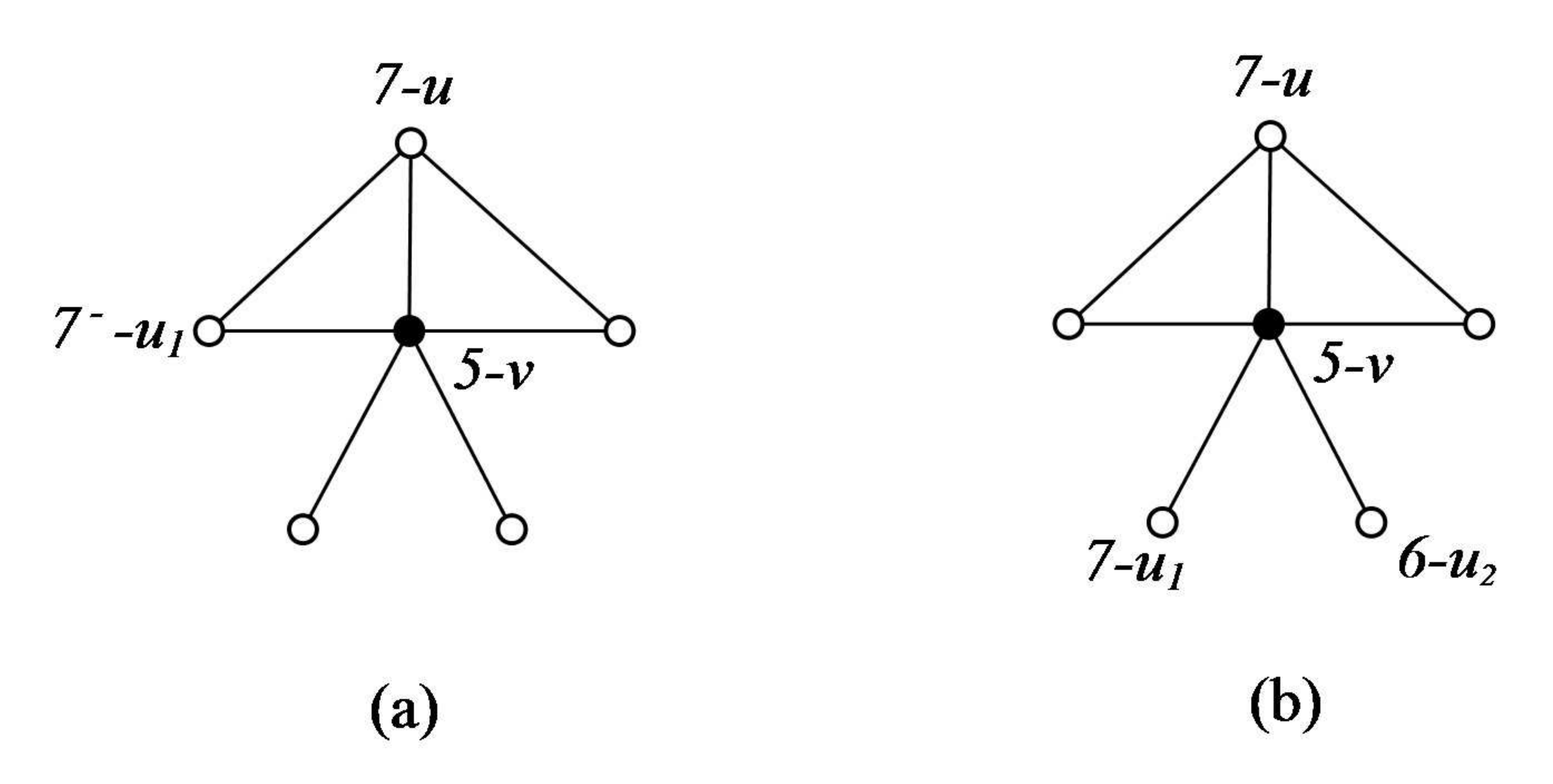}
\caption{The vertex $v$ is an $S_4$-neighbor of $u$.}
\label{S4}
\end{center}
\end{figure}

$\bullet$ The vertex $v$ is an $S_{4}$-$neighbor$ of $u$ with $d(u) = 7$ if $v$ is neither an $S_2$ nor $S_3$-neighbor of $u$, and one of the two following conditions is verified:

$\cdot$ $v$ is adjacent to a $7^-$-vertex $u_1$ such that $(u, v, u_1)$ are 3-faces (see Figure \ref{S4} (a)).

$\cdot$ $v$ is adjacent to a 7-vertex $u_1$ and a 6-vertex $u_2$ that are both distinct from $u$ (see Figure \ref{S4} (b)).

\section{Structural lemmas}
\begin{lemma}\label{l1}
Let $G$ be a planar graph with maximum degree $\Delta\geq8$ and $Y$ $(1 \leq|Y|\leq3)$ be a subset of nonadjacent vertices of $G$ on the same face $f_{0}$ such that $H= G- Y$ has at least one edge. 
Suppose that
\begin{description}
    \item[$(a)$] $d_{G}(v)\geq3$ for each vertex $v\in V(H)$, and
    \item[$(b)$] $d_{G}(u)+d_{G}(v)\geq\Delta+3$ for each edge $uv\in E(H)$.
\end{description}
Then $H$ contains one of the following configurations (see Figure \ref{lemma1}):
\begin{description}
  \item[$(1)$] there is a cycle $(u,v,w,x)$ such that $d_{G}(u)=d_{G}(w)=3$, $u,w\in V(H)$ and $v,x\in V(G)$;
  \item[$(2)$] $\Delta=8$,
  \begin{description}
  \item[$(2.1)$] there is a vertex $v\in V(H)$ with $d_{G}(v) = 8$ having three neighbors $u_{1}$, $u_{2}$, and $u_{3}$ in $H$ such that $u_{1}$ and $u_{2}$ are weak neighbors of $v$ with $d_{G}(u_{1}) = d_{G}(u_{2}) = 3$ and $d_{G}(u_{3})\leq5$;
  \item[$(2.2)$] there is a vertex $v\in V(H)$ with $d_{G}(v) = 8$ having four neighbors $u_{1}$, $u_{2}$, $u_{3}$ and $u_{4}$ in $H$ such that $u_{1}$ is a weak neighbor of $v$ and $u_{2}$ is a semiweak neighbor of $v$, with $d_{G}(u_{1}) = d_{G}(u_{2}) = 3$, $d_{G}(u_{3})\leq5$ and $d_{G}(u_{4})\leq5$;
  \item[$(2.3)$] there is a vertex $v\in V(H)$ with $d_{G}(v) = 8$ having four weak neighbors $u_{1}$, $u_{2}$, $u_{3}$ and $u_{4}$ in $H$ with $d_{G}(u_{1}) =3, d_{G}(u_{2})=d_{G}(u_{3})=4$ and $d_{G}(u_{4})\leq5$;
  \item[$(2.4)$] there is a vertex $v\in V(H)$ with $d_{G}(v) = 8$ having five neighbors $u_{1}$, $u_{2}$, $u_{3}$, $u_{4}$ and $u_{5}$ in $H$ such that $u_{1}$ is a weak neighbor of $v$ with $d_{G}(u_{1}) =3, d_{G}(u_{2}) = 4$, $d_{G}(u_{3})\leq5$, $d_{G}(u_{4})\leq5$ and $d_{G}(u_{5})\leq7$;
  \item[$(2.5)$] there is a vertex $v\in V(H)$ with $d_{G}(v) = 8$ having four weak neighbors $u_{1}$, $u_{2}$, $u_{3}$ and $u_{4}$ in $H$ such that $d_{G}(u_{1}) =3$, the vertex $u_{2}$ is an $E_{2}$-neighbor of $v$, $d_{G}(u_{3})\leq5$, and $d_{G}(u_{4})\leq5$;
  \item[$(2.6)$] there is a vertex $v\in V(H)$ with $d_{G}(v) = 7$ having three neighbors $u_{1}$, $u_{2}$ and $u_{3}$ in $H$ such that $u_{2}$ is adjacent to $u_{1}$ and $u_{3}$, $d_{G}(u_{2}) =6, d_{G}(u_{1}) =d_{G}(u_{3})=5$, and there is a vertex $w$ of degree 6, distinct from $u_{2}$, that is adjacent to $u_{3}$;
  \item[$(2.7)$] there is a vertex $v\in V(H)$ with $d_{G}(v) = 7$ having three weak neighbors $u_{1}$, $u_{2}$ and $u_{3}$ in $H$ such that $d_{G}(u_{1}) = d_{G}(u_{2}) = 4$ and either $u_{3}$ is an $S_{2}$, $S_{3}$, or $S_{4}$-neighbor, or $d_{G}(u_{3})=4$;
  \item[$(2.8)$] there is a vertex $v\in V(H)$ with $d_{G}(v) = 7$ having three neighbors $u_{1}$, $u_{2}$ and $u_{3}$ in $H$ such that $d_{G}(u_{1}) =4$, the vertex $u_{2}$ is an $S_{3}$-neighbor of $v$, $d_{G}(u_{3})\leq5$;
  \item[$(2.9)$] there is a vertex $v\in V(H)$ with $d_{G}(v) = 5$ having three neighbors $u_{1}$, $u_{2}$ and $u_{3}$ in $H$ such that $u_{2}$ is adjacent to $u_{1}$ and $u_{3}$ and $d_{G}(u_{1})= d_{G}(u_{2}) =d_{G}(u_{3})=6$.
  \end{description}
\end{description}

\begin{figure}[htbp]
\begin{center}
\includegraphics[scale=0.2]{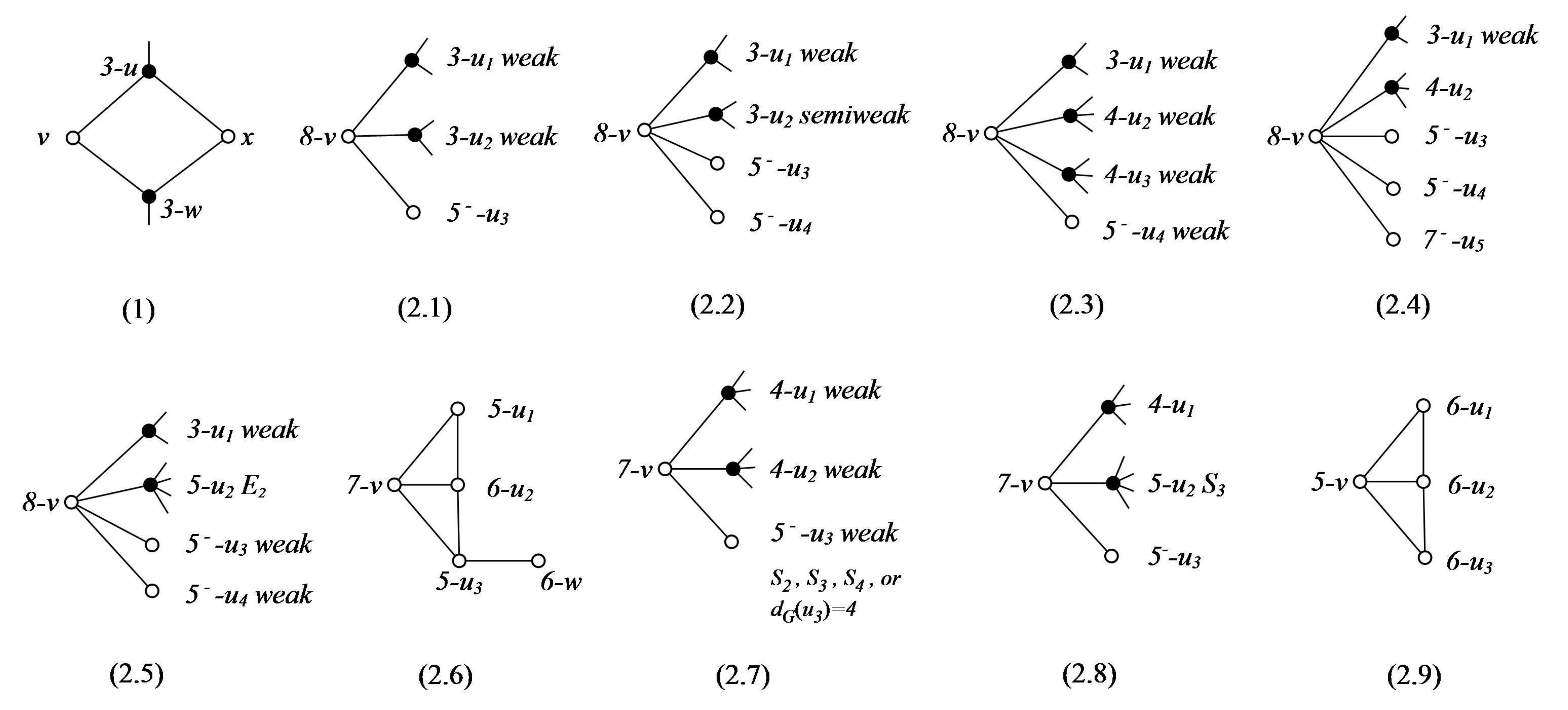}
\caption{unavoidable configurations of Lemma \ref{l1}.}
\label{lemma1}
\end{center}
\end{figure}

\end{lemma}

\begin{proof}
The proof is carried out by contradiction. Let $G$ be a counterexample to the lemma. By Euler's formula $|V|-|E|+|F|\geq 2$, we have
 $$\sum\limits_{v\in V(G)}(d(v)-6)+\sum\limits_{f\in F(G)}(2d(f)-6)\leq -12. $$
The above inequality can be written as
 $$\sum\limits_{v\in V(G)}(d(v)-6)+\sum\limits_{f\in F(G)\backslash f_0}(2d(f)-6)+(2d(f_0)+6)\leq 0.$$

We define $ch$ to be the {\em initial charge} by letting

\begin{equation*}
ch(x)=
\begin{cases}
d(x)-6,& x\in V(G),\\
2d(x)-6,& x\in F(G)\backslash f_0, \\
2d(x)+6, & x=f_0.
\end{cases}
\end{equation*}

For simplicity, let $VF(G)=V(G)\cup F(G)$. Thus, we have

$$\sum_{x\in VF(G)}ch(x)\leq 0.$$

Let $v\in V(G)$, $u\in V(H)$ and $f$ be a face of $G$. Now we define the discharging rules as follows.

{\it
\begin{description}
  \item[$\bf R1.$] Let $f=f_0$ and $v$ be incident with $f_0$. If $v\in Y$, then $f_0$ sends $6$ to $v$; otherwise let $Z$ be the set of vertices adjacent to $v$ and incident with $f_0$. If $|Z\backslash Y|=1$, then $f_0$ sends $1$ to $v$. If $|Z\backslash Y|=2$, then $f_0$ sends $2$ to $v$.
  \item[$\bf R2.$] Let $f\neq f_0$ and $u$ be incident with $f_0$ and $d_G(u)\leq5$. If $d(f)=4$, then $f$ sends $1$ to $u$; otherwise if $d(f)\geq5$, then $f$ sends $2$ to $u$.
  \item[$\bf R3.$] Let $v\in Y$. If $u$ is adjacent to $v$, then $v$ sends $1$ to $u$.
  \item[$\bf R4.$] Let $v\in V(H)$ and $d_G(v)\geq7$. If $v$ has a weak neighbor $u$ such that $d_G(u)=3$, then $v$ sends $1$ to $u$. If $v$ has a semiweak neighbor $u$ such that $d_G(u)=3$, then $v$ sends $\frac{1}{2}$ to $u$. If $v$ has a weak neighbor $u$ such that $d_G(u)=4$, then $v$ sends $\frac{1}{2}$ to $u$.
  \item[$\bf R5.$] Let $v\in V(H)$ and $d_G(v)\geq9$. If $v$ has a weak neighbor $u$ such that $d(u)=5$, then $v$ sends $\frac{1}{2}$ to $u$.
  \item[$\bf R6.$] Let $v\in V(H)$ and $d_G(v)=8$. If $v$ has an $E_{2}$-neighbor $u$, then $v$ sends $\frac{1}{2}$ to $u$. If $v$ has an $E_{3}$-neighbor $u$, then $v$ sends $\frac{1}{3}$ to $u$. If $v$ has an $E_{4}$-neighbor $u$, then $v$ sends $\frac{1}{4}$ to $u$.
  \item[$\bf R7.$] Let $v\in V(H)$ and $d_G(v)=7$. If $v$ has an $S_{2}$-neighbor $u$, then $v$ sends $\frac{1}{2}$ to $u$. If $v$ has an $S_{3}$-neighbor $u$, then $v$ sends $\frac{1}{3}$ to $u$. If $v$ has an $S_{4}$-neighbor $u$, then $v$ sends $\frac{1}{4}$ to $u$.
\end{description}}

Let $ch'(x)$ be the new charge according to the above discharging rules for each $x\in VF(G)$.  Since our rules only move charges around and do not affect the sum, we have

$$\sum_{x\in VF(G)}ch'(x)=\sum_{x\in VF(G)}ch(x)\leq 0.$$

In the following,  we shall show that $ch'(x)\geq 0$ for each $x\in VF(G)$ and $\sum_{x\in VF(G)}ch'(x)>0$ to obtain a contradiction.

\vspace{3mm}
Let $f$ be a face of $G$.

\vspace{3mm}
\noindent
\textbf{Case 1.} $f=f_0$.

\vspace{3mm}
R1 is equivalent to that for any $y\in Y$, there is a vertex incident with $f_0$ receives nothing from $f_0$. So $ch'(f_0) =ch(f_0)-|Y|\times6-(d(f_0)-2|Y|)\times2\geq 0$.

\vspace{3mm}
\noindent
\textbf{Case 2.} $f\not=f_0$.

\vspace{3mm}
By $(b)$, $f$ is incident with at most $\lfloor \frac{d(f)}{2} \rfloor$ $5^{-}$-vertices.

\romannumeral1. $d(f)= 3$. Then $ch'(f)=ch(f)=2d(f)-6=0$

\romannumeral2. $d(f)= 4$. Then $f$ is incident with at most two $5^{-}$-vertices. So $ch'(f)\geq ch(f)-2\times1=0$ by R2.

\romannumeral3. $d(f)= 5$. Then $f$ is incident with at most two $5^{-}$-vertices. So $ch'(f)\geq ch(f)-2\times2=0$ by R2.

\romannumeral4. $d(f)\geq 6$. Then $f$ is incident with at most $\lfloor \frac{d(f)}{2} \rfloor$ $5^{-}$-vertices. So $ch'(f)\geq ch(f)-\lfloor \frac{d(f)}{2} \rfloor \cdot 2 \geq d(f)-6\geq0$ by R2.

\vspace{3mm}
We have obtained that $ch'(f)\geq0$ for each $f\in F(G)$. Then we show that $ch'(v)\geq0$ for each $v\in V(G)$. If $v\in Y$, then $ch'(v)\leq ch(v)+6-d(v)=0$ by R1 and R3. So in the following, assume that $v\in V(H)$. We define that a vertex $u$ is a $Y_{f_0}$ $neighbor$ of $v$ when $u\in V(H)$ and one of the two faces adjacent to the edge $uv$ is $f_0$. Then we consider the following cases.


\vspace{3mm}
\noindent
\textbf{Case 1.} $d_G(v)=3$.

\vspace{3mm}
By $(b)$, the neighbors of $v$ are $\Delta$-vertices or in $Y$. Suppose that there are three faces $f_{1}$, $f_{2}$, and $f_{3}$ incident with $v$ and $d(f_{1})\geq d(f_{2})\geq d(f_{3})$. Let $N_G(v)=\{v_1, v_2, v_3\}$ such that $vv_{i}$ and $vv_{i+1}$ are incident with the face $f_i$ for any $1\leq i \leq2$ and $vv_{3}$ and $vv_{1}$ are incident with the face $f_3$. Note that if $f_0\in \{f_{1}, f_{2}, f_{3}\}$, R1 and R3 are equivalent to that $v$ receives 2 from $f_0$ and the $Y_{f_0}$ neighbors of $v$. So in the following, it can be considered $v$ receives $2$ from $f_0$.

\romannumeral1. $d(f_{1})\geq 5$ and  $d(f_{2})\geq 4$. If $f_0\in \{f_{1}, f_{2}, f_{3}\}$, then $v$ receives at least 2 from $f_0$ and at least 1 from $f_1$ or $f_2$; otherwise, $v$ receives 2 from $f_1$ and at least 1 from $f_2$. So $ch'(v)\geq ch(v)+2+1=0$ by R1, R2 and R3.

\romannumeral2. $d(f_{1})= d(f_{2})= d(f_{3})=4$. If $f_0\in \{f_{1}, f_{2}, f_{3}\}$, then $v$ receives at least 2 from $f_0$ and at least 2 from the other two faces; otherwise, $v$ receives 1 from each $f_i$, $1\leq i \leq3$. So $ch'(v)\geq ch(v)+3=0$ by R1, R2 and R3.

\romannumeral3. $d(f_{1})= d(f_{2})=4$ and  $d(f_{3})=3$. If $f_0\in \{f_{1}, f_{2}, f_{3}\}$, then $v$ receives at least 2 from $f_0$ and at least 1 from the other two faces. So $ch'(v)\geq ch(v)+2+1=0$ by R1, R2 and R3. Otherwise, $v$ receives 1 from each $f_i$, $1\leq i \leq2$. Besides, if there is a vertex $v_i \in Y$, $1\leq i \leq3$, then $v$ receives 1 from it by R3; otherwise $v$ is a semiweak neighbor of $v_1$ and $v_3$, then $v$ receives $\frac{1}{2}\times2$ from $v_1$ and $v_3$ by R4. So $ch'(v)\geq ch(v)+2+1=0$.

\romannumeral4. $d(f_{1})\geq 5$ and $d(f_{2})=d(f_{3})=3$. If $f_0\in \{f_{2}, f_{3}\}$, then $v$ receives at least 2 from $f_0$ and at least 2 from $f_1$. So $ch'(v)\geq ch(v)+2+2>0$ by R1, R2 and R3. Otherwise $v$ receives at least 2 from $f_1$. Besides, $v_2 \in V(H)$ and $v$ is a weak neighbor of $v_2$ or $v_2 \in Y$, then $v$ receives 1 from $v_2$ by R3 and R4. So $ch'(v)\geq ch(v)+2+1=0$.

\romannumeral5. $d(f_{1})=4$ and $d(f_{2})=d(f_{3})=3$. If $f_0\in \{f_{2}, f_{3}\}$, then $v$ receives at least 2 from $f_0$ and at least 1 from $f_1$; otherwise if $f_0=f_{1}$, $v$ receives at least 2 from $f_0$ and 1 from $v_2$ by R3 and R4. So $ch'(v)\geq ch(v)+2+1=0$. Suppose that $f_0$ is not incident with $v$. Then $v$ receives at least 1 from $f_1$. Besides, for any $i$, $1\leq i \leq3$, if $v_i \in Y$, then $v$ receives 1 from it; otherwise $v$ is a weak neighbor of $v_2$ and a semiweak neighbor of $v_1$ and $v_3$. In either case, $v$ receives at least 2 from the vertices in $N_G(v)$. So $ch'(v)\geq ch(v)+1+2=0$.

\romannumeral6. $d(f_{i})=3$, $1\leq i \leq3$. If $v_i \in Y$, then $v$ receives 1 from $v_i$; otherwise $v$ is a weak neighbor of $v_i$, $1\leq i \leq3$. In either case, $v$ receives at least 3 from the vertices in $N_G(v)$. So $ch'(v)\geq ch(v)+3=0$.

\vspace{3mm}
\noindent
\textbf{Case 2.} $d_G(v)=4$.

\vspace{3mm}
By $(b)$, the neighbors of $v$ are $(\Delta-1)^{+}$-vertices or in $Y$. Suppose that there are four faces $f_{1}$, $f_{2}$, $f_{3}$, and $f_{4}$ incident with $v$. Let $N_G(v)=\{v_1, v_2, v_3, v_4\}$.

Note that if $f_0\in \{f_{1}, f_{2}, f_{3}, f_{4}\}$, it is similar to case 1 that $v$ receives 2 from $f_0$. We assume $f_0=f_1$. So if there is a $4^+$-face $f_i$, $2\leq i\leq4$, then $v$ receives at least 1 from $f_i$ by R2; otherwise the faces except $f_1$ are all 3-face, assume that $vv_1$ and $vv_2$ are incident with $f_0$, then $v$ receives at least $\frac{1}{2}\times2$ from $v_2$ and $v_3$ by R4, so $ch'(v)\geq ch(v)+2+1>0$.

In the following, assume that $f_0$ is not incident with $v$ and $d(f_{1})\geq d(f_{2})\geq d(f_{3})\geq d(f_{4})$.

\romannumeral1. $d(f_{1})\geq d(f_{2})\geq 4$. The vertex $v$ is incident with at most two triangles, then $v$ receives at least 2 from $f_1$ and $f_2$ by R2. So $ch'(v)\geq ch(v)+2=0$.

\romannumeral2. $d(f_{1})\geq 4$ and $d(f_{i})=3$, $2\leq i \leq4$. The vertex $v$ is incident with exactly three triangles, then $v$ receives at least 1 from $f_1$ by R2. Besides, if there is a vertex $v_i \in Y$, $1\leq i \leq4$, then $v$ receives 1 from $v_i$ by R3; otherwise assume that $vv_1$ and $vv_2$ are incident with $f_1$, then $v$ is a weak neighbor of $v_3$ and $v_4$, thus $v$ receives $\frac{1}{2}\times2$ from $v_3$ and $v_4$ by R4. So $ch'(v)\geq ch(v)+1+1=0$.

\romannumeral3. $d(f_{i})=3$, $1\leq i \leq4$. Then the vertex $v$ is incident with four triangles. For any $i$, if $v_i \in Y$, then $v$ receives 1 from $v_i$; otherwise $v$ is a weak neighbor of $v_i$, $1\leq i \leq4$. In either case, $v$ receives at least 4 from the vertices in $N_G(v)$. So $ch'(v)\geq ch(v)+4>0$.

\vspace{3mm}
\noindent
\textbf{Case 3.} $d_G(v)=5$.

\vspace{3mm}
By $(b)$, the neighbors of $v$ are $(\Delta-2)^{+}$-vertices or in $Y$. Suppose that there are five faces $f_{1}$, $f_{2}$, $f_{3}$, $f_{4}$ and $f_{5}$ incident with $v$. Let $N_G(v)=\{v_i, 1\leq i\leq5\}$ and such that $vv_{i}$ and $vv_{i+1}$ are incident with the face $f_i$ for any $1\leq i \leq4$ and $vv_{5}$ and $vv_{1}$ are incident with the face $f_5$.

Note that if there is an vertex $v_i \in Y$, $1\leq i \leq5$, then $v$ receives at least 1 from $v_i$ by R3, so $ch'(v)\geq ch(v)+1=0$; otherwise if $f_0\in \{f_{1}, f_{2}, f_{3}, f_{4}, f_5\}$, then $v$ receives at least 2 from $f_0$ by R1, so $ch'(v)\geq ch(v)+2>0$; otherwise if there is a face $f_i$ such that $d(f_i)\geq4$, $1\leq i \leq5$, then $v$ receives at least 1 from $f_i$ by R2, so $ch'(v)\geq ch(v)+1=0$.

In the following, assume that $v_i \notin Y$, $f_0$ is not incident with $v$ and $d(f_{i})=3$, $1\leq i \leq5$. So $v$ is incident with five triangles.

\vspace{1mm}
\textbf{\uppercase\expandafter{\romannumeral1}.} $\Delta\geq9$. By $(b)$, $d(v_{i})\geq7$, $1\leq i \leq5$.

\vspace{1mm}
\romannumeral1. The vertex $v$ is adjacent to at least two $9^{+}$-vertices.

Then $v$ receives $\frac{1}{2}\times2$ from them by R5. So $ch'(v)\geq ch(v)+\frac{1}{2}\times2=0$.

\romannumeral2. The vertex $v$ is adjacent to exactly one $9^{+}$-vertex.

Without loss of generality, we assume $d(v_{1})\geq 9$. Then $v$ receives $\frac{1}{2}$ from $v_1$ by R5.

(\romannumeral2-1) The vertex $v$ is adjacent to four 8-vertices. Then $v$ is an $E_4$-neighbor of $v_{3}$ and $v_{4}$, so $v$ receives $\frac{1}{4}\times2$ from $v_3$ and $v_{4}$ by R6. Thus $ch'(v)\geq ch(v)+\frac{1}{2}+\frac{1}{4}\times2=0$.

(\romannumeral2-2) The vertex $v$ is adjacent to three 8-vertices and one 7-vertices. We assume w.l.o.g. that $d(v_{2})= d(v_{3})=8$. If $d(v_{4})=8$ and $d(v_{5})=7$, then $v$ is an $E_4$-neighbor of $v_{3}$ and an $E_3$-neighbor of $v_{4}$, so $v$ receives $\frac{1}{4}$ from $v_3$ and $\frac{1}{3}$ from $v_{4}$ by R6; otherwise $d(v_{4})=7$ and $d(v_{5})=8$, then $v$ is an $E_3$-neighbor of $v_{3}$ and $v_{5}$, so $v$ receives $\frac{1}{3}\times2$ from $v_3$ and $v_{5}$ by R6. Thus $ch'(v)\geq ch(v)+\frac{1}{2}+\frac{1}{3}+ \frac{1}{4}>0$.

(\romannumeral2-3) The vertex $v$ is adjacent to two 8-vertices and two 7-vertices. If $d(v_{2})= d(v_{3})=8$ and $d(v_{4})=d(v_{5})=7$, then $v$ is an $E_3$-neighbor of $v_{3}$ and an $S_4$-neighbor of $v_{4}$, so $v$ receives $\frac{1}{3}$ from $v_3$ and $\frac{1}{4}$ from $v_{4}$ by R6 and R7. If $d(v_{2})= d(v_{4})=8$ and $d(v_{3})=d(v_{5})=7$, then $v$ is an $E_3$-neighbor of $v_{2}$ and $v_{4}$, so $v$ receives $\frac{1}{3}\times2$ from $v_2$ and  $v_{4}$ by R6. If $d(v_{2})= d(v_{5})=8$ and $d(v_{3})=d(v_{4})=7$, then $v$ is an $E_3$-neighbor of $v_{2}$ and $v_{5}$,  and an $S_4$-neighbor of $v_{3}$ and $v_{4}$, so $v$ receives $\frac{1}{3}\times2$ from $v_{2}$ and $v_{5}$, and $\frac{1}{4}\times2$ from $v_{3}$ and $v_{4}$ by R6 and R7. If $d(v_{3})= d(v_{4})=8$ and $d(v_{2})=d(v_{5})=7$, then $v$ is an $E_3$-neighbor of $v_{3}$ and $v_{4}$, so $v$ receives $\frac{1}{3}\times2$ from $v_3$ and $v_{4}$ by R6. Thus $ch'(v)\geq ch(v)+\frac{1}{2}+\frac{1}{3}+ \frac{1}{4}>0$.

(\romannumeral2-4) The vertex $v$ is adjacent to one 8-vertices and three 7-vertices. We assume w.l.o.g. that $d(v_{2})= d(v_{3})=7$. If $d(v_{4})=8$ and $d(v_{5})=7$, then $v$ is an $E_3$-neighbor of $v_{4}$ and an $S_4$-neighbor of $v_{2}$ and $v_3$, so $v$ receives $\frac{1}{3}$ from $v_4$ and $\frac{1}{4}\times2$ from $v_{2}$ and $v_3$ by R6 and R7; otherwise $d(v_{4})=7$ and $d(v_{5})=8$, then $v$ is an $E_3$-neighbor of $v_{5}$ and an $S_4$-neighbor of $v_{2}$, $v_3$ and $v_4$, so $v$ receives $\frac{1}{3}$ from $v_5$ and $\frac{1}{4}\times3$ from $v_{2}$,$v_3$ and $v_4$ by R6 and R7. Thus $ch'(v)\geq ch(v)+\frac{1}{2}+\frac{1}{3}+ \frac{1}{4}\times2>0$.

(\romannumeral2-5) The vertex $v$ is adjacent to four 7-vertices. Then $v$ is an $S_4$-neighbor of $v_{i}$, $2\leq i \leq5$, so $v$ receives $\frac{1}{4}\times4$ from them by R7. Thus $ch'(v)\geq ch(v)+\frac{1}{2}+\frac{1}{4}\times4>0$.

\romannumeral3. The vertex $v$ is adjacent to no $9^{+}$-vertex.

(\romannumeral3-1) The vertex $v$ is adjacent to at least four 8-vertices. We assume w.l.o.g. that $d(v_{i})=8$, $1\leq i \leq4$. Then $v$ is an $E_3$ or $E_4$-neighbor of $v_{i}$, so $v$ receives at least $\frac{1}{4}$ from $v_i$, $1\leq i \leq4$ by R6. Thus $ch'(v)\geq ch(v)+\frac{1}{4}\times4=0$.

(\romannumeral3-2) The vertex $v$ is adjacent to three 8-vertices and two 7-vertices. We assume w.l.o.g. that $d(v_1)=7$ and $d(v_{2})= d(v_{3})=8$. If $d(v_{4})=8$ and $d(v_{5})=7$, then $v$ is an $S_4$-neighbor of $v_{1}$ and $v_5$, an $E_3$-neighbor of $v_2$ and $v_{4}$, and an $E_4$-neighbor of $v_3$. So $v$ receives $\frac{1}{4}\times3$ from $v_{1}$, $v_3$ and $v_5$, $\frac{1}{3}\times2$ from $v_2$ and $v_{4}$ by R6 and R7; otherwise $d(v_{4})=7$ and $d(v_{5})=8$, then $v$ is an $E_3$-neighbor of $v_2$, $v_{3}$ and $v_{5}$. So $v$ receives $\frac{1}{3}\times3$ from them by R6. Thus $ch'(v)\geq ch(v)+\frac{1}{3}\times3=0$.

(\romannumeral3-3) The vertex $v$ is adjacent to two 8-vertices and three 7-vertices. We assume w.l.o.g. that $d(v_1)=8$ and $d(v_{2})= d(v_{3})=7$. If $d(v_{4})=7$ and $d(v_{5})=8$, then $v$ is an $S_4$-neighbor of $v_{2}$, $v_3$ and $v_4$, and an $E_3$-neighbor of $v_1$ and $v_{5}$. So $v$ receives $\frac{1}{4}\times3$ from $v_{2}$, $v_3$ and $v_4$, and $\frac{1}{3}\times2$ from $v_1$ and $v_{5}$ by R6 and R7; otherwise $d(v_{4})=8$ and $d(v_{5})=7$, then $v$ is an $S_4$-neighbor of $v_2$ and $v_{3}$, and an $E_3$-neighbor of $v_1$ and $v_4$. So $v$ receives $\frac{1}{4}\times2$ from $v_{2}$ and $v_3$, and $\frac{1}{3}\times2$ from $v_1$ and $v_{4}$ by R6 and R7. Thus $ch'(v)\geq ch(v)+\frac{1}{4}\times2+\frac{1}{3}\times3=0$.

(\romannumeral3-4) The vertex $v$ is adjacent to at least four 7-vertices. We assume w.l.o.g. that $d(v_{i})=7$, $1\leq i \leq4$. For each $i$, $v$ is an $S_4$-neighbor of $v_{i}$, so $v$ receives $\frac{1}{4}$ from $v_i$ by R6. Thus $ch'(v)\geq ch(v)+\frac{1}{4}\times4=0$

\vspace{1mm}
\textbf{\uppercase\expandafter{\romannumeral2}.} $\Delta=8$. By $(b)$, $d(v_{i})\geq6$, $1\leq i \leq5$.

\vspace{1mm}
\romannumeral1. The vertex $v$ is adjacent to at least three 6-vertices.

By (2.9), they are exactly three. Without loss of generality, we assume $d(v_{1})= d(v_{2})= d(v_{4})=6$, hence $d(v_{3}), d(v_{5})\geq7$. Then $v$ is an $S_2$ or $E_2$ -neighbor of $v_{3}$ and $v_{5}$, so $v$ receives $\frac{1}{2}\times2$ from them by R6 and R7. Thus $ch'(v)\geq ch(v)+\frac{1}{2}\times2=0$.

\romannumeral2. The vertex $v$ is adjacent to exactly two 6-vertices.

(\romannumeral2-1) The vertex $v$ is adjacent to two consecutive 6-vertices. Without loss of generality, we assume $d(v_{1})= d(v_{2})=6$ and $d(v_{3})\geq d(v_{5})$. If $d(v_{3})= d(v_{5})=8$, then $v$ is an $E_2$-neighbor of $v_{3}$ and $v_{5}$, so $v$ receives $\frac{1}{2}\times2$ from $v_3$ and $v_{5}$ by R6. If $d(v_{3})=8$ and $d(v_{5})=7$, then $v$ is an $E_2$-neighbor of $v_{3}$, an $S_3$ or $S_4$-neighbor of $v_{5}$, and an $S_4$ or $E_3$-neighbor of $v_{4}$, so $v$ receives $\frac{1}{2}$ from $v_3$ and at least $\frac{1}{4}\times2$ from $v_4$ and $v_5$ by R6 and R7. If $d(v_{3})= d(v_{5})=7$, then $v$ is an $S_3$-neighbor of $v_{3}$ and $v_{5}$,  and an $S_3$ or $E_3$-neighbor of $v_{4}$, so $v$ receives $\frac{1}{3}\times3$ from $v_{3}$, $v_4$ and $v_{5}$. Thus $ch'(v)\geq ch(v)+1=0$.

(\romannumeral2-2) The vertex $v$ is not adjacent to two consecutive 6-vertices. We assume w.l.o.g. that $d(v_{1})= d(v_{4})=6$ and $d(v_{2})\geq d(v_{3})$. If $d(v_{3})=8$, then $v$ is an $E_3$ or $S_2$-neighbor of $v_{5}$ and an $E_3$-neighbor of $v_{2}$ and $v_3$, so $v$ receives at least $\frac{1}{3}\times3$ from $v_2$, $v_3$ and $v_{5}$ by R6 and R7; otherwise $d(v_{3})=7$, then $v$ is an $E_2$ or $S_2$-neighbor of $v_{5}$ and an $S_3$, $S_4$, or $E_3$-neighbor of $v_{2}$ and $v_3$, so $v$ receives $\frac{1}{2}$ from $v_5$ and at least $\frac{1}{4}\times2$ from $v_2$ and $v_3$ by R6 and R7. Thus $ch'(v)\geq ch(v)+1=0$.

\romannumeral3. The vertex $v$ is adjacent to exactly one 6-vertex.

Without loss of generality, we assume $d(v_{1})=6$ and $d(v_{2})\geq d(v_{5})$ or  $d(v_{3})\geq d(v_{4})$ if  $d(v_{2})=d(v_{5})$. If $d(v_{5})=8$ and $d(v_{3})=d(v_{4})$, then $v$ is an $E_3$-neighbor of $v_{2}$ and $v_{5}$ and an $S_4$ or $E_4$-neighbor of $v_{3}$ and $v_{4}$, so $v$ receives $\frac{1}{3}\times2$ from $v_2$ and $v_{5}$ and $\frac{1}{4}\times2$ from $v_3$ and $v_{4}$ by R6. If $d(v_{5})=8$ and $d(v_{3})\neq d(v_{4})$, then $v$ is an $E_3$-neighbor of $v_{2}$, $v_{3}$ and $v_{5}$, so $v$ receives $\frac{1}{3}\times3$ from them by R6 and R7. If $d(v_{5})=7$, then $v$ is an $E_3$, $E_4$ or $S_4$-neighbor of $v_{i}$, $2\leq i\leq5$, so $v$ receives at least $\frac{1}{4}\times4$ from them by R6 and R7. Thus $ch'(v)\geq ch(v)+1=0$.

\romannumeral4. The vertex $v$ is adjacent to no 6-vertex.

(\romannumeral4-1) The vertex $v$ is adjacent to at least four 8-vertices. Then $v$ is an $E_3$ or $E_4$-neighbor of each of them, so $v$ receives at least $\frac{1}{4}\times4$ from them by R6. Thus $ch'(v)\geq ch(v)+1=0$.

(\romannumeral4-2) The vertex $v$ is adjacent to two consecutive 7-vertices. We assume w.l.o.g. that $d(v_{1})= d(v_{2})=7$. Then $v$ is an $S_4$-neighbor of $v_{1}$ and $v_2$, and an $S_4$ or $E_3$-neighbor of $v_{3}$ and $v_5$, so $v$ receives at least $\frac{1}{4}\times4$ from them by R6 and R7. Thus $ch'(v)\geq ch(v)+1=0$.

(\romannumeral4-3) The vertex $v$ is adjacent to at most three 8-vertices and no two consecutive 7-vertices. We assume w.l.o.g. that $d(v_{1})= d(v_{3})=7$ and  $d(v_{2})= d(v_{4})= d(v_{5})=8$ . Then $v$ is an $E_3$-neighbor of $v_{2}$, $v_4$ and $v_5$, so $v$ receives $\frac{1}{3}\times3$ from them by R6 and R7. Thus $ch'(v)\geq ch(v)+1=0$.

\vspace{3mm}
\noindent
\textbf{Case 4.} $d_G(v)=6$.

\vspace{3mm}
The vertex $v$ gives nothing away, so $ch'(v)\geq ch(v=0$. Note that if $f_0$ is incident with $v$, R1 and R3 are equivalent to that $v$ receives 2 from $f_0$ and the $Y_{f_0}$ neighbors of $v$, then $ch'(v)\geq ch(v)+2>0$.

\vspace{3mm}
\noindent
\textbf{Case 5.} $d_G(v)=7$.

\vspace{3mm}
By $(b)$, the neighbors of $v$ are $4^{+}$-vertices or in $Y$ and $v$ has at most three weak neighbors of degree at most 5. Suppose that there are seven faces $f_{i}$, $1\leq i\leq7$, incident with $v$. Let $N_G(v)=\{v_i, 1\leq i\leq7\}$ and such that $vv_{i}$ and $vv_{i+1}$ are incident with the face $f_i$ for any $1\leq i \leq6$ and $vv_{7}$ and $vv_{1}$ are incident with the face $f_7$.

Note that if there is a vertex $v_i \in Y$, $1\leq i \leq7$, then $v$ receives at least 1 from $v_i$ by R3 and gives at most $\frac{1}{2}\times3$ to other vertices by R4 and R7, so $ch'(v)\geq ch(v)+1-\frac{1}{2}\times3>0$; otherwise if $f_0$ is incident with $v$, then $v$ receives at least 2 from $f_0$ by R1 and gives at most $\frac{1}{2}\times3$ to other vertices by R4 and R7, so $ch'(v)\geq ch(v)+2-\frac{1}{2}\times3>0$.

In the following, assume that $v_i \notin Y$ and $f_0$ is not incident with $v$.

\vspace{1mm}
\textbf{\uppercase\expandafter{\romannumeral1}.} $\Delta\geq9$. By $(b)$, $d(v_{i})\geq5$, $1\leq i \leq7$.

\vspace{1mm}
The vertex $v$ has at most three weak neighbors of degree 5 and they are not $S_2$ or $S_3$-neighbor of $v$ by $(b)$. Then $v$ gives at most $\frac{1}{4}\times3$ to them by R7, so $ch'(v)\geq ch(v)-\frac{1}{4}\times3>0$.

\vspace{1mm}
\textbf{\uppercase\expandafter{\romannumeral2}.} $\Delta=8$. By $(b)$, $d(v_{i})\geq4$, $1\leq i \leq7$.

\vspace{1mm}
\romannumeral1. The vertex $v$ has an $S_2$-neighbor $v_1$.

So $d(v_{2})=d(v_{7})=6$ by definition of an $S_2$-neighbor. By (2.6), if $v_3$ or $v_6$ is a weak neighbor of $v$, then $d(v_3)>5$ or $d(v_6)>5$ ,respectively. We assume w.l.o.g. that $d(v_{4})\geq d(v_{5})$. If $v_4$ and $v_5$ are adjacent, then $d(v_{4})>5$ by $(b)$. Thus $v$ has at most two weak neighbors of degree at most 5 ($v_1$ and possibly $v_5$) and $v$ gives at most $\frac{1}{2}\times2$ to them by R4 and R7. So $ch'(v)\geq ch(v)-\frac{1}{2}\times2=0$.

\romannumeral2. The vertex $v$ has at least two weak neighbors of degree 4.

By (2.7), $v$ is adjacent to no other weak neighbor of degree 4 and no $S_2$, $S_3$, or $S_4$-neighbor. Thus $v$ gives at most $\frac{1}{2}\times2$ to them by R4. So $ch'(v)\geq ch(v)-\frac{1}{2}\times2=0$.

\romannumeral3. The vertex $v$ has exactly one weak neighbor $v_1$ of degree 4 and no $S_2$-neighbor.

If $v$ has an $S_3$-neighbor $v_i$, $2\leq i\leq7$, then other neighbors are all $6^+$-vertices by (2.8). Thus $v$ gives $\frac{1}{2}$ to $v_1$ and $\frac{1}{3}$ to $v_i$ by R4 and R7. So $ch'(v)\geq ch(v)-\frac{1}{2}+\frac{1}{3}>0$. Otherwise if $v$ has no $S_3$-neighbor, then $v$ has at most two other weak neighbors of degree 5 by assumption. Thus $v$ gives $\frac{1}{4}\times2$ to them and $\frac{1}{2}$ to $v_1$ by R4 and R7. So $ch'(v)\geq ch(v)-\frac{1}{4}\times2+\frac{1}{2}=0$.

\romannumeral4. The vertex $v$ has no weak neighbor of degree 4 and no $S_2$-neighbor.

Then $v$ has at most three weak neighbors of degree 5 by assumption. Thus $v$ gives at most $\frac{1}{3}\times3$ to them by R7. So $ch'(v)\geq ch(v)-\frac{1}{3}\times3=0$.

\vspace{3mm}
\noindent
\textbf{Case 6.} $d_G(v)=8$.

\vspace{3mm}
By $(b)$, the neighbors of $v$ are $3^{+}$-vertices or in $Y$ and $v$ has at most four semiweak neighbors of degree 3 or weak neighbors of degree at most 5. 
Let $N_G(v)=\{v_i, 1\leq i\leq8\}$ 

\vspace{1mm}
\textbf{\uppercase\expandafter{\romannumeral1}.} $\Delta\geq9$. By $(b)$, $d(v_{i})\geq4$, $1\leq i \leq8$.

\vspace{1mm}
The vertex $v$ has at most four weak neighbors of degree 4 or 5. Then $v$ gives at most $\frac{1}{2}\times4$ to them by R4 and R6, so $ch'(v)\geq ch(v)-\frac{1}{2}\times4=0$. Note that if $f_0$ is incident with $v$, R1 and R3 are equivalent to that $v$ receives at least 2 from $f_0$ and the $Y_{f_0}$ neighbors of $v$, then $ch'(v)\geq ch(v)+2-\frac{1}{2}\times4>0$

\vspace{1mm}
\textbf{\uppercase\expandafter{\romannumeral2}.} $\Delta=8$.

\vspace{1mm}
Note that if there is exactly a vertex $v_i \in Y$ and $f_0$ is not incident with $v$, $1\leq i \leq8$, then $v$ has either three weak neighbors of degree 3, or two weak neighbors of degree 3 and two semiweak neighbors of degree 3 or weak neighbors of degree 4 or 5 by (1). Thus $v$ receives 1 from $v_i$ by R3 and gives at most 3 to other vertices by R4 and R6, so $ch'(v)\geq ch(v)+1-3=0$. Otherwise if there is at least two vertices $v_i \in Y$, $1\leq i \leq8$, or $f_0$ is incident with $v$, then $v$ has at most three weak neighbors of degree 3. Thus $v$ receives at least 2 by R1 and R3, and gives at most $1\times3+\frac{1}{2}$ to other vertices by R4 and R6, so $ch'(v)\geq ch(v)+2-1\times3-\frac{1}{2}>0$.

In the following, assume that $v_i \notin Y$ and $f_0$ is not incident with $v$.

\romannumeral1. The vertex $v$ has at least two weak neighbors $v_1$ and $v_2$ of degree 3.

By (2.1), $v$ has exactly two neighbors of degree at most 5 ($v_1$ and $v_2$) and $v$ gives $1\times2$ to them by R4. So $ch'(v)\geq ch(v)-1\times2=0$.

\romannumeral2. The vertex $v$ has exactly one weak neighbor $v_1$ of degree 3 and at least one semiweak neighbor $v_2$ of degree 3.

By (2.2), $v$ is adjacent to at most one other neighbor $v_3$ of degree at most 5 and $v_3$ is not a weak neighbor of degree 3. Thus $v$ gives at most $\frac{1}{2}\times2+1$ to $v_1$, $v_2$ and $v_3$ by R4 and R6. So $ch'(v)\geq ch(v)-\frac{1}{2}\times2-1=0$.

\romannumeral3. The vertex $v$ has exactly one weak neighbor $v_1$ of degree 3, no semiweak neighbor of degree 3, and at least two weak neighbor $v_2$ and $v_3$ of degree 4.

By (2.3), $v$ has no other weak neighbor of degree at most 5, then $v$ gives $1$ to $v_1$, $\frac{1}{2}\times2$ to  $v_2$ and $v_3$ by R4 and R6. So $ch'(v)\geq ch(v)-1-\frac{1}{2}\times2=0$.

\romannumeral4. The vertex $v$ has exactly one weak neighbor $v_1$ of degree 3, no semiweak neighbor of degree 3, exactly one weak neighbor $v_2$ of degree 4, and at least one $E_2$ or $E_3$-neighbor $v_3$ .

The vertices $v$ and $v_3$ have a common neighbor $v_4 $ of degree 6 or 7 by (b). Then by (2.4), $v$ has no other neighbor of degree at most 5. Thus $v$ gives $1$ to $v_1$, $\frac{1}{2}$ to $v_2$ and at most $\frac{1}{2}$ to $v_3$ by R4 and R6. So $ch'(v)\geq ch(v)-1-\frac{1}{2}-\frac{1}{2}=0$.

\romannumeral5. The vertex $v$ has exactly one weak neighbor $v_1$ of degree 3, no semiweak neighbor of degree 3, exactly one weak neighbor $v_2$ of degree 4, and no $E_2$ or $E_3$-neighbor.

The vertices $v$ has at most two other weak neighbor  $v_3$ and $v_4 $ of degree at most 5. Then $v$ gives $1$ to $v_1$, $\frac{1}{2}$ to $v_2$, and $\frac{1}{4}\times2$ to $v_3$ and $v_4$ by R4 and R6. So $ch'(v)\geq ch(v)-1-\frac{1}{2}-\frac{1}{4}\times2=0$.

\romannumeral6. The vertex $v$ has exactly one weak neighbor $v_1$ of degree 3, no semiweak neighbor of degree 3, no weak neighbor of degree 4, and at least one $E_2$-neighbor $v_2$.

By (2.5), $v$ has at most one other weak neighbor $v_3$ of degree 5, then $v$ gives $1$ to $v_1$ and at most $\frac{1}{2}\times2$ to $v_2$ and $v_3$ by R4 and R6. So $ch'(v)\geq ch(v)-1-\frac{1}{2}\times2=0$.

\romannumeral7. The vertex $v$ has exactly one weak neighbor $v_1$ of degree 3, no semiweak neighbor of degree 3, no weak neighbor of degree 4, and no $E_2$-neighbor.

The vertex $v$ has at most three other weak neighbors $v_2$, $v_3$ and $v_4$ of degree 5. Since $v$ has no $E_2$-neighbor, they are  $E_3$ or $E_4$-neighbor of $v$. Thus $v$ gives $1$ to $v_1$ and at most $\frac{1}{3}\times3$ to $v_2$ and $v_3$ by R4 and R6. So $ch'(v)\geq ch(v)-1-\frac{1}{3}\times3=0$.

\romannumeral8. The vertex $v$ has no weak neighbor of degree 3.

The vertex $v$ has at most four neighbors of degree at most 5. Thus $v$ gives at most $\frac{1}{2}\times4$ to them by R4 and R6. So $ch'(v)\geq ch(v)-\frac{1}{2}\times4=0$.

\vspace{3mm}
\noindent
\textbf{Case 7.} $d_G(v)\geq9$.

\vspace{3mm}
By $(b)$, The vertex $v$ has at most $\lfloor\frac{1}{2}d_G(v)\rfloor$ semiweak neighbors with degree 3 or weak neighbors of degree at most 5. 
Let $N_G(v)=\{v_i, 1\leq i\leq9\}$ 

Note that if there is at least two vertices $v_i \in Y$, $1\leq i \leq9$ or $f_0$ is incident with $v$, then $v$ receives at least 2 by R1 and R3, and gives at most $1\times\lfloor\frac{1}{2}d_G(v)\rfloor$ to other vertices by R4 and R5, so $ch'(v)\geq ch(v)+2-1\times\lfloor\frac{1}{2}d_G(v)\rfloor>0$. Otherwise if there is exactly one  vertex $v_i \in Y$, $1\leq i \leq9$, then $ch'(v)\geq ch(v)+2-1\times\lfloor\frac{1}{2}d_G(v)\rfloor\geq0$.
In the following, assume that $v_i \notin Y$ and $f_0$ is not incident with $v$.

\romannumeral1. The vertex $v$ has no weak neighbor of degree 3.

The vertex $v$ has at most $\lfloor\frac{1}{2}d_G(v)\rfloor$ neighbors of degree at most 5. Thus $v$ gives at most $\frac{1}{2}\times\lfloor\frac{1}{2}d_G(v)\rfloor$ to them by R4 and R5. So $ch'(v)\geq ch(v)-\frac{1}{2}\times\lfloor\frac{1}{2}d_G(v)\rfloor>0$.

\romannumeral2. The vertex $v$ has at least one weak neighbor of degree 3.

By (b), $d_G(v)=\Delta$.

(\romannumeral2-1) $\Delta=9$. Assume that $v$ is adjacent to $k$ weak neighbors of degree 3 and $k\leq3$ by (b) and (1). If $k=3$, then $v$ has no other neighbor of degree at most 5. Thus $v$ gives at most $1\times3$ to them by R4. So $ch'(v)\geq ch(v)-1\times3=0$. Otherwise if $k\leq2$, then $v$ has at most $(4-k)$ other neighbors of degree at most 5. Thus $v$ gives at most $1\times k+\frac{1}{2}\times(4-k)$ to them by R4 and R5. So $ch'(v)\geq ch(v)-(1\times k+\frac{1}{2}\times(4-k))\geq0$.

(\romannumeral2-2) $\Delta=10$. Assume that $v$ is adjacent to $k$ weak neighbors of degree 3 and $k\leq3$ by (b) and (1). Then $v$ has at most $(5-k)$ other neighbors of degree at most 5. Thus $v$ gives at most $1\times k+\frac{1}{2}\times(5-k)$ to them by R4 and R5. So $ch'(v)\geq ch(v)-(1\times k+\frac{1}{2}\times(5-k))\geq0$.

(\romannumeral2-3) $\Delta\geq11$. The vertex $v$ has at most $\lfloor\frac{1}{2}d_G(v)\rfloor$ neighbors of degree at most 5. Thus $v$ gives at most $1\times\lfloor\frac{1}{2}d_G(v)\rfloor$ to them by R4 and R5. So $ch'(v)\geq ch(v)-1\times\lfloor\frac{1}{2}d_G(v)\rfloor\geq0$.

\vspace{3mm}
Till now, we have checked that $ch'(x)\geq 0$ for any element $x\in VF(G)$. Now we begin to find a vertex or a face $x\in VF(G)$ such that $ch'(x) > 0$. If $|Y|\leq2$, then $ch'(f_0)>0$. So we assume $|Y|=3$. Let $v\in V(H)$ be a vertex incident with $f_0$. According to Case 2 - Case 7, $ch'(v)>0$ if $d_G(v)\geq4$. So we only need to consider the Case 1 that $d_G(v)=3$.

We assume w.l.o.g. that $f_0=f_1$. R1 and R3 are equivalent to that $v$ receives at
least 2 from $f_0$ and the $Y_{f_0}$ neighbors of $v$. If $f_2$ or $f_3$ is a $5^+$-face, then $ch'(v)\geq ch(v)+2+2>0$ by R2. If $f_2$ and $f_3$ are both 4-face, then $ch'(v)\geq ch(v)+2+1\times2>0$ by R2. If $f_2$ and $f_3$ is a 4-face and a 3-face, then $v$ is a semiweak neighbor of $v_2$ or $v_2\in Y$, so $ch'(v)\geq ch(v)+2+1+\frac{1}{2}>0$ by R2, R3 and R4. So we suppose that they are 3-faces. Then $v$ is a weak neighbor of $v_2$ or $v_2\in Y$, so $v$ receives 1 from $v_2$ by R3 and R4. In addition, if $v_1 \notin Y$ or $v_3 \notin Y$, then $v$ is a semiweak neighbor of it, so  $ch'(v)\geq ch(v)+2+1+\frac{1}{2}>0$ by R4. So we consider the case that $v_1, v_3 \in Y$, $f_0=f_1$, and $f_2$ and $f_3$ are 3-faces. Then $d_G(v_2)=\Delta$ by $(b)$ and $v_2$ receives 2 from $v_1$ and $v_3$, so it follows from Case 6 and Case 7 that $ch'(v_2)>0$.

Hence we complete the proof of the lemma.
\end{proof}

\begin{figure}[htbp]
\begin{center}
\includegraphics[scale=0.4]{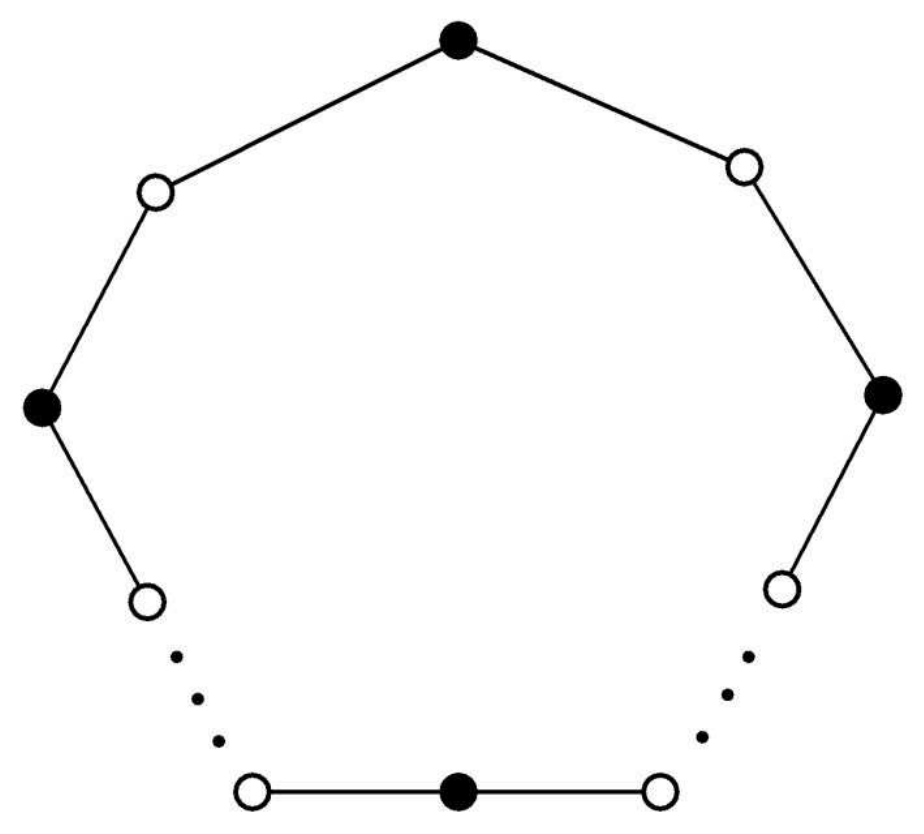}
\caption{A $2$-$alternating$ $cycle$}
\label{2alter}
\end{center}
\end{figure}

We define a $2$-$alternating$ $cycle$ in a graph $G$ which is a cycle of even length in which alternate vertices have degree 2 in $G$ (see Figure \ref{2alter}). Accordingly, we define a $3$-$alternator$ to be a bipartite subgraph $F$ of $G$ with partite sets $U$, $W$ such that, for each $u\in U$, $2\leq d_{F}(u)= d_{G}(u)\leq3$, and for each $w\in W$, either $d_{F} (w)\geq3$ or $w$ has exactly two neighbours in $U$, both with degree exactly $14-d_{G}(w)$ (this last being possible only if $d_{G}(w)=11$ or 12).

\begin{lemma}\label{l2}
Let $G$ be a planar graph of maximum degree $\Delta(G)\geq12$ and $Y$ $(1 \leq|Y|\leq3)$ be a subset of nonadjacent vertices of $G$ on the same face $f_{0}$ such that $H= G-Y$ has at least one edge.
Suppose that
\begin{description}
    \item[$(a)$] $d_{G}(v)\geq2$ for each vertex $v\in V(H)$, and
    \item[$(b)$] $H$ contains no $2$-alternating cycle or $3$-alternator.
\end{description}
Then there is an edge $uv\in E(H)$ such that $d_{G}(u)+d_{G}(v)\leq\Delta(G)+1$.
\end{lemma}

\begin{proof}

The proof is carried out by contradiction. Let $G$ be a counterexample to the lemma with as few vertices as possible and, subject to this, as many edges as possible. So for each edge $uv\in E(H)$, $d_{G}(u)+d_{G}(v)\geq14$.

Let $V_{f_{0}}$ be a set of the vertices that are incident with $f_{0}$ and $G^{*}$ be the graph obtained by deleting all 2-vertices from $G\setminus V_{f_{0}}$ along with their incident edges.

\begin{claim}
Except $f_0$ and the faces $(a,u,b,v)$ such that $a, b \in Y$, the faces of $G^{*}$ are $3$-faces.
\end{claim}

\begin{proof}
The neighbours $u$, $w$ of a 2-vertex $v$ in $G\setminus V_{f_{0}}$ have degree at least 12 or in $Y$, and so $v$ is incident with a 3-face in $G$, otherwise we could join $u$, $w$ by an edge, which is a contradiction with the maximality of edges. In addition to $f_0$ and the faces $(a,u,b,v)$ that $a, b \in Y$, we suppose that there is a $4^{+}$-face $f$ in $G^{*}$. Let $x$, $y$, $z$ be three consecutive vertices in the boundary of $f$, chosen so that $|\{x,z\}\cap Y|\leq1$, $y\notin Y$ and $d_{G}(y)$ is as small as possible. If $d_{G}( y)\leq 6$, then $d_G(x)\geq8$ or $x\in Y$, and $d_G(z)\geq8$ or $z\in Y$; otherwise, $d_G(x)\geq d_G( y)\geq7$ or $x\in Y$, and $d_G(z)\geq7$ or $z\in Y$. In either case we could add the edge $xz$ inside $f$ without violating any of the hypotheses of the theorem, and this contradiction proves the claim.
\end{proof}

Let $U_{G}=\{v\in V(H), d_{G}(v)\leq3\}$ and $U_{Y}=\{v\in V(H), d_{G}(v)\leq3$ and $N_G(v)\subseteq Y\}$. So we consider $U=U_{G}\setminus U_{Y}$ and $W=N_H(U)$. If $X\subseteq E(H)$, write $U(X)$, $W(X)$ for the vertices in $U$, $W$ that are incident with the edges in $X$. Let $X^{\ast}$ be a maximal subset of the edges between $U$ and $W$ such that
\begin{description}
    \item[$(1)$] for each $u\in U(X^{\ast})$, $d_{X^{\ast}}(u)=1$,
    \item[$(2)$] for each $w\in W(X^{\ast})$, $d_{X^{\ast}}(w)\leq2$ and $w$ has at most one $X^{\ast}$-neighbour
in $U$ of degree $14-d_{G}(w)$, and
    \item[$(3)$] $N(U\setminus U(X^{\ast}))\subseteq W\setminus W(X^{\ast})$.
\end{description}

Since $H$ contains no 3-alternator, if $U\setminus U(X^{\ast})\neq\emptyset$, then there is a vertex $w\in N(U\setminus U(X^{\ast}))$ such that if we add to $X^{\ast}$ the edges between $U\setminus U(X^{\ast})$ and $w$, then the resulting set satisfies (1) to (3), a contradiction. Thus $U(X^{\ast})=U$. If $uw\in X$ ($u\in U, w\in W$) then we call $w$ the $master$ of $u$ and $u$ a $dependent$ of $w$.

 By Euler's formula $|V|-|E|+|F|\geq 2$, we have
 $$\sum\limits_{v\in V(G)}(d(v)-6)+\sum\limits_{f\in F(G)}(2d(f)-6)\leq -12, $$
which can be written as follows
 $$\sum\limits_{v\in V(G)}(d(v)-6)+\sum\limits_{f\in F(G)\backslash f_0}(2d(f)-6)+(2d(f_0)+6)\leq 0.$$

We define $ch$ to be the {\em initial charge} by letting

\begin{equation*}
ch(x)=
\begin{cases}
d(x)-6,& x\in V(G),\\
2d(x)-6,& x\in F(G)\backslash f_0, \\
2d(x)+6, & x=f_0.
\end{cases}
\end{equation*}

Let $VF(G)=V(G)\cup F(G)$. Thus, we have

$$\sum_{x\in VF(G)}ch(x)\leq 0.$$

Note that each 2-vertex is incident with at most one 3-face. Let $v\in V(G)$, $u\in V(H)$ and $f$ be a face of $G$. Now we define the discharging rules as follows.

{\it
\begin{description}
  \item[$\bf R1.$] Let $f=f_0$ and $v$ be incident with $f_0$. If $v\in Y$, then $f_0$ sends $6$ to $v$; otherwise let $Z$ be the set of vertices adjacent to $v$ and incident with $f_0$. If $|Z\backslash Y|=1$, then $f_0$ sends $1$ to $v$. If $|Z\backslash Y|=2$, then $f_0$ sends $2$ to $x$.
  \item[$\bf R2.$] Let $v\in Y$. If $u$ is adjacent to $v$, then $v$ sends $1$ to $u$.
  \item[$\bf R3.$] Let $d_G(u)=2$, then $u$ receives $2$ from its master vertex and $2$ from the $4^{+}$-face except $f_0$ which is incident with it.
  \item[$\bf R4.$] Let $d_G(u)=3$, then $u$ receives $2$ from its master vertex and $\frac{1}{2}$ from each of its neighbours in $H$.
  \item[$\bf R5.$] Let $4\leq d_G(u)\leq5$, then $u$ receives $\frac{1}{2}$ from each of its neighbours in $H$.
\end{description}}

Let $ch'(x)$ be the new charge according to the above discharging rules for each $x\in VF(G)$.  Since our rules only move charges around and do not affect the sum, we have

$$\sum_{x\in VF(G)}ch'(x)=\sum_{x\in VF(G)}ch(x)\leq 0.$$

In the following,  we shall show that $ch'(x)\geq 0$ for each $x\in VF(G)$ and $\sum_{x\in VF(G)}ch'(x)>0$ to obtain a contradiction.

\vspace{3mm}
Let $f$ be a face of $G$.

\vspace{3mm}
\noindent
\textbf{Case 1.} $f=f_0$.

\vspace{3mm}
R1 is equivalent to that for any $y\in Y$, there is a vertex incident with $f_0$ receives nothing from $f_0$. So $ch'(f_0) =ch(f_0)-|Y|\times6-(d(f_0)-2|Y|)\times2\geq 0$.

\vspace{3mm}
\noindent
\textbf{Case 2.} $f\not=f_0$.

\vspace{3mm}
\romannumeral1. $d(f)= 3$. Then $ch'(f)=ch(f)=2d(f)-6=0$

\romannumeral2. $d(f)= 4$. Then $f$ is incident with at most one 2-vertices. So $ch'(f)\geq ch(f)-2=0$ by R3.

\romannumeral3. $d(f)= 5$. Then $f$ is incident with at most two 2-vertices. So $ch'(f)\geq ch(f)-2\times2=0$ by R3.

\romannumeral4. $d(f)\geq 6$. Then $f$ is incident with at most $\lfloor \frac{d(f)}{2} \rfloor$ 2-vertices. So $ch'(f)\geq ch(f)-\lfloor \frac{d(f)}{2} \rfloor \times 2 \geq d(f)-6\geq0$ by R3.

\vspace{3mm}
We have obtained that $ch'(f)\geq0$ for each $f\in F(G)$. Then we show that $ch'(v)\geq0$ for each $v\in V(G)$. If $v\in Y$, then $ch'(v)\leq ch(v)+6-d(v)=0$ by R1 and R2. So in the following, we assume that $v\in V(H)$.
If $v$ has given $\frac{1}{2}$ to a neighbour $u$, then we call $u$ a \emph{half-neighbour} of $v$. It is clearly that $d(u)\leq5$ and  $d(v)\geq9$.

\vspace{3mm}
\noindent
\textbf{Case 1.} $d_G(v)=2$.

\vspace{3mm}
$N_G(v)\subseteq Y$, then $v$ receives 2 from $Y$ by R2. Besides, if $v$ is incident with $f_0$, then it receives 2 from the other $4^{+}$-face by R3, so $ch'(v)\geq ch(v)+2+2=0$; otherwise $v$ receives $2\times2$ from two $4^{+}$-faces by R3, so $ch'(v)\geq ch(v)+2+2\times2>0$. If $|N_G(v)\cap Y|=1$, then $v$ receives 1 from $Y$, 2 from its master vertex by R2 and at least 1 from the faces  which is incident with it by R1 and R3, so $ch'(v)\geq ch(v)+1+2+1=0$. Otherwise $v$ receives 2 from its master vertex and at least 2 from $4^{+}$-faces which is incident with it by R3. So $ch'(v)\geq ch(v)+2+2=0$.

\vspace{3mm}
\noindent
\textbf{Case 2.} $d_G(v)=3$.

\vspace{3mm}
If $N_G(v)\subseteq Y$, then $v$ receives 3 from $Y$ by R2. Otherwise $v$ receives 2 from its master vertex and at least $\frac{1}{2}\times2$ from the other two neighbors by R4. So $ch'(v)\geq ch(v)+3=0$.

\vspace{3mm}
\noindent
\textbf{Case 3.} $d_G(v)=4$.

\vspace{3mm}
Then $v$ receives at least $\frac{1}{2}\times4$ from its neighbors by R5. So $ch'(v)\geq ch(v)+2=0$.

\vspace{3mm}
\noindent
\textbf{Case 4.} $d_G(v)=5$.

\vspace{3mm}
Then $v$ receives at least $\frac{1}{2}\times5$ from its neighbors by R5. So $ch'(v)\geq ch(v)+2.5>0$.

\vspace{3mm}
\noindent
\textbf{Case 5.} $d_G(v)=6$.

\vspace{3mm}
Then $ch'(v)\geq ch(v)=0$.

\vspace{3mm}
\noindent
\textbf{Case 6.} $d_G(v)=7$ or 8.

\vspace{3mm}
Then $ch'(v)\geq ch(v)>0$.

\vspace{3mm}
\noindent
\textbf{Case 7.} $d_G(v)=9$ or 10.

\vspace{3mm}
Then $v$ gives at most $\frac{1}{2}\times5$ to half-neighbours. So $ch'(v)\geq ch(v)-\frac{1}{2}\times5>0$.

\vspace{3mm}
\noindent
\textbf{Case 8.} $d_G(v)=11$.

\vspace{3mm}
Then $v$ gives at most 2 to a dependent of degree 3 and at most $\frac{1}{2}\times5$ to half-neighbours. So $ch'(v)\geq ch(v)-2-\frac{1}{2}\times5>0$.

\vspace{3mm}
\noindent
\textbf{Case 9.} $d_G(v)=12$.

\vspace{3mm}
Then $v$ gives at most $2\times2$ to dependents, but at most one of these is 2-vertex, and at most $\frac{1}{2}\times4$ to half-neighbours. So $ch'(v)\geq ch(v)-2\times2-\frac{1}{2}\times4=0$.

\vspace{3mm}
\noindent
\textbf{Case 10.} $d_G(v)\geq13$.

\vspace{3mm}
Then $v$ gives at most $2\times2$ to dependents and at most $\frac{1}{2}\times\lfloor\frac{d(v)-2}{2}\rfloor$ to half-neighbours. So $ch'(v)\geq ch(v)-2\times2-\frac{1}{2}\times\lfloor\frac{d(v)-2}{2}\rfloor>0$.

\vspace{3mm}
Till now, we have checked that $ch'(x)\geq 0$ for any element $x\in VF(G)$. Since $\sum_{x\in VF(G)}$
$ch'(x)\leq 0$, it follows that $ch'(x)=0$ and $d(v)\in\{2, 3, 4, 6, 12\}$ for each vertex $v\in V(H)$. Moreover, each 12-vertex has two dependents, four half-neighbours, six neighbours of degree 12 (by the claim), no neighbors in $Y$ and it is not incident with $f_0$ (otherwise it receives at least 1 by R1 and R2). Besides, each 2-vertex is adjacent to two 12-vertices or two vertices in $Y$. It follows that there are no 2-vertices that have two 12-neighbors except $Y$, since no 12-vertex is adjacent to a nondependent 2-vertex. Moreover, the vertices in $\{v\in V(H)$, $d_G(v)=12\}$ induce a 6-regular subgraph of G, which must be a triangulation by Euler's formula. Hence there are no 4-vertices or 6-vertices either, and every 12-vertex is adjacent to six 3-vertices, which is contradicted to that there is no 3-alternator.

Hence we complete the proof of the lemma.
\end{proof}

\section{Some structural properties of $K_5$-minor free graphs}
Let $G$ be a connected graph, $T$ be a tree, and $ \mathcal{F}=\{V_t\subset V(G):t\in T\}$ be a family of subsets of $V(G)$. The ordered set $(T, \mathcal{F})$ is called a $tree$-$decomposition$ of $G$ if it satisfies the following conditions:

(T1) $V(G)=\bigcup_{t\in T}V_{t}$;

(T2) for each edge $e\in E(G)$, there exists a vertex $t\in V(T)$ such that the two end-vertices of $e$ are included in $V_t$;

(T3) if $t_1,t_2,t_3\in V(T)$ and $t_2$ is on the $(t_1,t_3)$-path of $T$, then $V_{t_1}\cap V_{t_3}\subset V_{t_2}$.

Give a fixed tree-decomposition $(T, \mathcal{F})$ of $G$, for each pair of adjacent vertices $s$ and $t$ in $T$, $V_s\cap V_t$ form a vertex cut of $G$, which is called a $separator$ $set$ of $(T, \mathcal{F})$. The graph $G_t=G[V_t]$ $(t\in T)$ is called the $part$ of $(T, \mathcal{F})$.  If the induced subgraph of any separator set of $(T, \mathcal{F})$ is a complete graph, $(T, \mathcal{F})$ is called  $simple$. Moreover, if every separator set of a simple $(T, \mathcal{F})$ has at most $k$ vertices, then $(T, \mathcal{F})$ is called $k$-$simple$.

\begin{figure}[htbp]
\begin{center}
\includegraphics[scale=0.6]{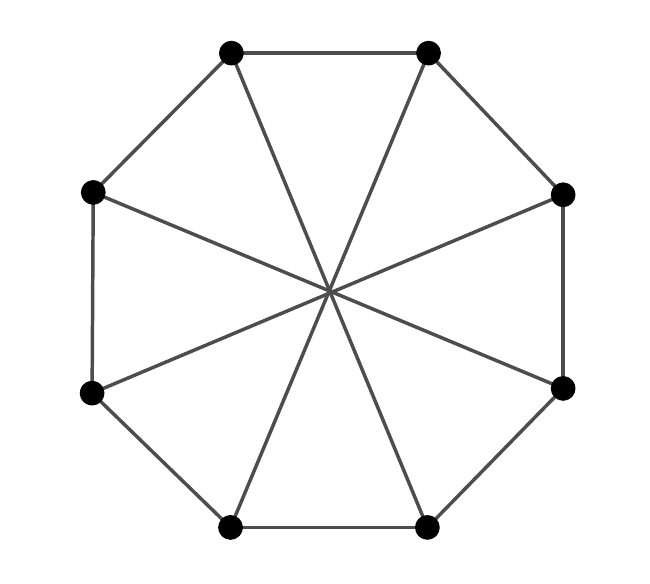} \\
\caption{The Wagner graph $W$}
\label{f}
\end{center}
\end{figure}

\begin{lemma}{\rm \cite{W}}\label{K5}
Let $G$ be an edge-maximal graph without a $K_5$ minor. If $|G|\geq 4$ then $G$ has a $3$-simple tree-decomposition $(T, \mathcal{F})$ such that each part is  a planar triangulation or the Wagner graph $W$$ $(see Figure $\ref{f})$.
\end{lemma}

The above lemma implies that  every $K_5$-minor free graph has a tree-decomposition $(T, \mathcal{F})$ such that each part is a planar graph or the Wagner graph and each separator set has size at most $3$.

\begin{lemma}\label{k5}
Let $G$ be a $K_5$-minor free graph and $\Delta(G)\geq8$. Then one of the following conditions holds:
\begin{description}
  \item[$(1)$] $\delta(G)\leq3$;
  \item[$(2)$] for each edge $uv\in E(G)$, $d_{G}(u)+d_{G}(v)\leq\Delta+2$;
  \item[$(3)$] $G$ has a subgraph isomorphic to one of configurations in Figure \ref{lemma1}.
\end{description}
\end{lemma}

\begin{proof}
Suppose to be contrary that $G$ is a counterexample to the lemma such that $|V(G)|$ is as small as possible. Then \begin{description}
  \item[$(a)$] $\delta(G)\geq2$;
  \item[$(b)$] for each edge $uv\in E(G)$, $d_{G}(u)+d_{G}(v)\geq\Delta+3$;
  \item[$(c)$] $G$ contains no configuration in Figure \ref{lemma1}.
\end{description}

Let  $(T, \mathcal{F})$  be a tree-decomposition of $G$ such that each part is a planar graph or the Wagner graph, each separator set is of size at most $3$ and $|V(T)|$ is as small as possible. Suppose that $|V(T)|=1$. Then $G$ must be a planar graph and $|V(G)|\geq 5$ (since the Wagner graph is of the maximum degree $3$, which is a contradiction to $(b)$). Let $v\in V(G)$ and $Y=\{v\}$. Then $G$ satisfies the conditions $(a)$ and $(b)$ of Lemma \ref{l1}, it follows that $G$ satisfies the lemma, a contradiction to $(c)$. So $|V(T)|\geq 2$.

Let $v_1v_2...v_t$ $(t\geq 2)$ be a longest path of $T$. Then $v_1$ is a leaf of $T$. By $(b)$, $G_{v_1}$ is a planar graph. Let $S_{12}=V_{v_1}\cap V_{v_2}$, $G'_1=G_{v_1}\backslash S_{12}$ and $G^*_1=G'_1\cup\{xy|\, x\in V(G'_1), y\in S_{12}$ and $xy\in E(G)\}$. Without loss of generality, we  assume that $G'_1$ is connected (for otherwise we can consider a connected component of $G'_1$). If $|V(G'_1)|\geq 2$, then it follows from Lemma \ref{l1} that $G$ contains a vertex satisfying the lemma, a contradiction. If $|V(G'_1)|=1$ and $|S_{12}|\leq 2$, then $\delta(G)\leq 2$, a contradiction to $(a)$. So $|V(G'_1)|=1$ and $|S_{12}|=3$, that is, $G^*_1$ is a star of order 4. Let $V(G^*_1)=\{u, u_{1}, u_{2}, u_{3}\}$ such that $\{u\}=V(G'_1)$ and $S_{12}=\{u_{1}, u_{2}, u_{3}\}$.  Then $d_G(u)=3$.

Since $|S_{12}|=3$, $G_{v_2}$ is a planar graph. Let $K=G_{v_2}\cup \{xy|\, x,y\in S_{12}\; \text{and}\; xy\notin G_{v_2}\}$. Then $K$ is also a planar graph.  We embed $K$ into the plane such that $S_{23}=V_{v_2}\cap V_{v_3}$ (if $t=2$, then we choose any vertex of $V_{v_2}$ as $S_{23}$) are located on the unbounded face $f_0$. By the minimality of $T$, The cycle $C=u_{1}u_{2}u_{3}u_{1}$ of $K$ must be a separated triangle and by the similar arguments as above, the inner part of $C$ is equivalent to a leaf of $T$ and is also a star of order $4$. Let $w$ be the inner vertex of $C$. Then $d_G(w)=3$ and $N_G(w)= N_{K}(w)= N_G(u)=\{u_{1}, u_{2}, u_{3}\}$. By (b), $d_G(u_{1})= d_G(u_{2})=d_G(u_{3})=\Delta$. Thus $G$ has the configuration (1) in Figure \ref{lemma1}, which is contradicted to $(c)$. We complete the proof of the lemma.
\end{proof}

\begin{lemma}\label{k52}
Let $G$ be a $K_5$-minor free graph and $\Delta(G)\geq12$. Then one of the following conditions holds:
\begin{description}
  \item[$(1)$] $\delta(G)<2$;
  \item[$(2)$] $G$ contains $2$-alternating cycle or $3$-alternator;
  \item[$(3)$] there is an edge $uv\in E(G)$, $d_{G}(u)+d_{G}(v)\leq\Delta(G)+1$.
\end{description}
\end{lemma}

\begin{proof}
Suppose to be contrary that $G_{1}$ is a counterexample to the lemma such that $|V(G_{1})|$ is as small as possible. Now we construct a new graph $G$ from $G_{1}$ by contracting all 2-vertices and deleting all
contracted multiple edges. Thus $G$ is also a $K_{5}$-minor free graph with $\delta(G_{1})\geq3$. Since the neighbors of any 2-vertex of $G_1$ have degree of $\Delta(G_1)$. Then we have
\begin{description}
  \item[$(a)$] $\delta(G)\geq3$;
  \item[$(b)$] $G$ contains no 2-alternating cycle or 3-alternator;
  \item[$(c)$] for any edge $uv\in E(G)$, $d_{G}(u)+d_{G}(v)\geq\Delta(G)+2$.
\end{description}

Let  $(T, \mathcal{F})$  be a tree-decomposition of $G$ such that each part is a planar graph or the Wagner graph, each separator set is of size at most $3$ and $|V(T)|$ is as small as possible. Suppose that $|V(T)|=1$. Then $G$ must be a planar graph and $|V(G)|\geq 5$ (since the Wagner graph is of the maximum degree $3$, it is contradicted to $(c)$). Let $v\in V(G)$ and $Y=\{v\}$. Then $G$ satisfies Lemma \ref{l2}, a contradiction. So $|V(T)|\geq 2$.

Let $v_1v_2...v_t$ $(t\geq 2)$ be a longest path of $T$. Then $v_1$ is a leaf of $T$. By $(c)$, $G_{v_1}$ is a planar graph. Let $S_{12}=V_{v_1}\cap V_{v_2}$, $G'_1=G_{v_1}\backslash S_{12}$ and $G^*_1=G'_1\cup\{xy|\, x\in V(G'_1), y\in S_{12}$ and $xy\in E(G)\}$. Without loss of generality, we assume that $G'_1$ is connected (for otherwise we can consider a connected component of $G'_1$). If $|V(G'_1)|\geq 2$, then it follows from Lemma \ref{l2} that $G$ satisfies the lemma, a contradiction. If $|V(G'_1)|=1$ and $|S_{12}|\leq 2$, then $\delta(G)\leq 2$, a contradiction to $(a)$. So $|V(G'_1)|=1$ and $|S_{12}|=3$, that is, $G^*_1$ is a star of order 4. Let $V(G^*_1)=\{u, u_{1}, u_{2}, u_{3}\}$ such that $\{u\}=V(G'_1)$ and $S_{12}=\{u_{1}, u_{2}, u_{3}\}$.  Then $d_G(u)=3$.

Since $|S_{12}|=3$, $G_{v_2}$ is a planar graph. Let $K=G_{v_2}\cup \{xy|\, x,y\in S_{12}\; \text{and}\; xy\notin G_{v_2}\}$. Then $K$ is also a planar graph.  We embed $K$ into the plane such that $S_{23}=V_{v_2}\cap V_{v_3}$ (if $t=2$, then we choose any vertex of $V_{v_2}$ as $S_{23}$) are located on the unbounded face $f_0$. By the minimality of $T$, The cycle $C=u_{1}u_{2}u_{3}u_{1}$ of $K$ must be a separated triangle and by the similar arguments as above, the inner part of $C$ is equivalent to a leaf of $T$ and is also a star of order $4$. Let $w$ be the inner vertex of $C$. Then $d_G(w)=3$ and $N_G(w)= N_{K}(w)= N_G(u)=\{u_{1}, u_{2}, u_{3}\}$. By(c), $d_G(u_{1})= d_G(u_{2})=d_G(u_{3})\geq\Delta-1$. Thus $G$ has the 3-alternator, which is contradicted to $(b)$. We complete the proof of the lemma.
\end{proof}

\section{The proof of Theorem \ref{th1} and Theorem \ref{th2}}

\begin{proof}[Proof of Theorem $\ref{th1}$] Suppose, to be contrary, that $G$ is a counterexample to
Theorem \ref{th1} with $|V| + |E|$ as small as possible. Then
\begin{description}
  \item[$(1)$] $\delta(G)\geq3$;
  \item[$(2)$] $d_{G}(u)+d_{G}(v)\geq\Delta+3$;
  \item[$(3)$] the configurations of Figure \ref{lemma1} are reducible, that is, they cannot be subgraphs of $G$.
\end{description}
The reducible configurations of (1)-(3) can be found in \cite{Bonamy} and \cite{Borodin}. By Lemma \ref{k5}, this is contradictory. We complete the proof of Theorem \ref{th1}.
\end{proof}

\begin{proof}[Proof of Theorem $\ref{th2}$] Suppose, to be contrary, that $G$ is a counterexample to
Theorem \ref{th2} with $|V| + |E|$ as small as possible. Then
\begin{description}
  \item[$(1)$] $\delta(G)\geq3$;
  \item[$(2)$] $G$ contains no 2-alternating cycle or 3-alternator;
  \item[$(3)$] for any edge $uv\in E(G)$, $d_{G}(u)+d_{G}(v)\geq\Delta(G)+2$.
\end{description}
The reducible configurations of (1)-(3) can be found in \cite{Borodin2}. By Lemma \ref{k52}, this is contradictory. We complete the proof of Theorem \ref{th2}.
\end{proof}

\end{document}